\documentclass{article}[12pt]
\usepackage{amsmath,amsfonts,amsthm}
\usepackage{amssymb}

\usepackage{epsfig}

\usepackage{subfigure}
\usepackage{graphics}
\usepackage{graphicx}
\usepackage{color}
\usepackage{soul}
\usepackage[makeroom]{cancel}

\usepackage{latexsym}
\usepackage{graphics}
\usepackage{amsmath}
\usepackage{amsthm}
\usepackage{xspace}
\usepackage{amssymb}
\usepackage{color}
\usepackage{cite}
\usepackage{latexsym}
\usepackage{graphics}
\usepackage{amsmath}
\usepackage{amsthm}
\usepackage{xspace}
\usepackage{amssymb}
\usepackage{epsfig}
\usepackage{comment}

\DeclareMathOperator*{\argmin}{arg\,min}

\newcommand{\beql}[1]{\begin{equation}\label{#1}}
\newcommand{\eeql}{\end{equation}}
\newcommand{\eqn}[1]{(\ref{#1})}

\newcommand{\R}{\mathbb{R}}
\newcommand{\pr}{\mathbb{P}}
\newcommand{\E}{\mathbb{E}}

\newcommand{\cx}{{\cal X}}
\newcommand{\cm}{{\cal M}}

\newcommand{\cc}{{\cal C}}

\newcommand{\bI}{{\bf I}}

\newtheorem{thm}{Theorem}
\newtheorem{lem}[thm]{Lemma}
\newtheorem{prop}[thm]{Proposition}

\newtheorem{definition}[thm]{Definition}
\newtheorem{conjecture}[thm]{Conjecture}

\theoremstyle{remark}

\begin{document}

\title{A large-scale particle system with independent jumps and distributed synchronization
}

\author
{
Yuliy Baryshnikov \\
Math and ECE Departments
and Coordinated Science Lab\\
University of Illinois at Urbana-Champaign\\
Urbana, IL 61801\\
\texttt{ymb@illinois.edu}\\
\and
Alexander L. Stolyar \\
ISE Department
and Coordinated Science Lab\\
University of Illinois at Urbana-Champaign\\
Urbana, IL 61801\\
\texttt{stolyar@illinois.edu}%\\
%\and 
%dhdhdhd
}

\date{\today}

\maketitle

\begin{abstract}

We study a system consisting of $n$ particles, moving forward in jumps on the real line.
Each particle can make both independent jumps, whose sizes have some distribution, or
``synchronization'' jumps, which allow it to join a randomly chosen other particle if the latter happens to be ahead of it.
System state is the empirical distribution of particle locations.
The mean-field asymptotic regime, where $n\to\infty$, is considered. 
We prove that $v_n$, the steady-state speed of the particle system advance, converges, as $n\to\infty$, to a limit $v_{**}$ which can be easily found from a {\em minimum speed selection principle.} Also, 
as $n\to\infty$, we prove the convergence of the system dynamics to that of a deterministic mean-field limit (MFL).
We show that the average speed of advance of any MFL is lower bounded by $v_{**}$,
and the speed of a ``benchmark'' MFL, resulting from all particles initially co-located, is equal to  $v_{**}$.

In the special case of exponentially distributed independent jump sizes, we prove that a traveling wave MFL with speed $v$ exists if and only if $v\ge v_{**}$, with $v_{**}$ having simple explicit form; we also show the existence of traveling waves for the modified systems, with a left or right boundary moving at a constant speed $v$. Using these traveling wave existence results, we provide bounds on an MFL average speed of advance, depending on the right tail exponent of its initial state. We conjecture that these results for exponential jump sizes generalize to general jump sizes.

\end{abstract}



%\newpage

\section{Introduction}

\subsection{Model, motivation and (informally) main results}

We consider a system consisting of $n$ particles, moving forward in jumps on the real line.
A particle can jump either independently or via ``synchronizations'' with other particles.
Independent jumps of a particle occur at time points of an independent Poisson process of constant rate $\lambda \ge 0$; the jump sizes are i.i.d., distributed as a random variable $Z>0$ with CDF $J(\cdot)$ and finite mean, $\E Z < \infty$;
without loss of generality we can and do assume $\E Z =1$. 
Each particle also has another independent Poisson process (of synchronizations), of constant rate $\mu>0$, at which points it chooses another particle uniformly at random, and makes a synchronization jump to the location of the chosen particle if it is ahead -- to the right -- of its own current location.
The system state at a given time is the empirical distribution of particle locations.
We focus on the mean-field asymptotic regime, where $n$ becomes large.

Our model is of the type, where particles move both autonomously (independent jumps) and according to a synchronization mechanism. (See \cite{balazs-racz-toth-2014, GMP97, GSS96, jonck-fpp-2020, Mal-Man-2006, Malyshkin-2006, Manita-2006, Manita-2009, Manita-2014, Man-Sch-2005} and references therein for previous models of this type.) Such models are motivated by distributed systems, where agents need to both evolve and synchronize
their states, and the synchronization is done in distributed fashion, via random peer-to-peer communication. 
The special feature of our model is that, even in the asymptotic limit, the autonomous movement of a particle is discontinuous (consists of random jumps). This brings up the question of how much
the synchronization mechanism affects the overall system speed of advance, and what is a good trade-off
between “self-propulsion” (independent jumps) rate and synchronization rate. In models of technological systems,
such as parallel computing/simulation or wireless systems (cf.  \cite{GMP97, GSS96, Mal-Man-2006, Malyshkin-2006, Manita-2006, Manita-2009, Manita-2014, Man-Sch-2005}), 
particles' locations correspond to agents' local ``clocks,'' which advance locally, but need to be synchronized. Our specific system can also model a business environment, where multiple companies can improve their products either via own
investments (“self-propulsion”) or via acquisitions (“synchronization”).

The main question that we address in this paper is very typical for this kind of systems: as $n$ becomes large, what is the average speed at which the system state (empirical distribution) moves forward. Specifically, we will be interested in two metrics of the speed of advance: the steady-state speed of advance $v_n$, especially $\lim_{n\to\infty} v_n$;
the average speed of the {\em mean-field limit} (MFL), which, informally speaking, the deterministic limit of the system state dynamics as $n\to\infty$. Our main results are as follows.

\begin{itemize}

\item We prove that $v_n \to v_{**}$ as ${n\to\infty}$, where the value of $v_{**}$ can be easily found from a {\em minimum speed selection principle.} 

\item As $n\to\infty$, we prove the convergence of the system dynamics to that of a deterministic MFL. We also prove this for two modified systems, with a left or right boundary moving at a constant speed $v$. 

\item We show that the average speed of advance of any MFL is lower bounded by $v_{**}$,
and the speed of a ``benchmark'' MFL (BMFL), resulting from all particles initially co-located, is equal to  $v_{**}$.

\item In the special case of exponentially distributed independent jump sizes:

\begin{itemize}

\item We prove that a traveling wave MFL with speed $v$ exists if and only if $v\ge v_{**}$, with $v_{**}$ having very simple explicit form. In addition, we show that a traveling wave with moving left boundary exists for any boundary speed $v > v_{**}$, and a traveling wave with moving right boundary exists for any boundary speed $v < v_{**}$. 

\item Using monotonicity and the traveling waves as bounds, we obtain bounds on an MFL average speed of advance, depending on the right tail exponent of the initial state.

\end{itemize}

\end{itemize}

We conjecture that the above results for exponential jump sizes generalize to general jump sizes as well.

Using our results, in a large-scale system, one can easily optimize a trade-off between $\lambda$ and $\mu$
(between ``self-propulsion'' and synchronization efforts) to achieve maximum speed. For example, one can maximize $v_{**} = v_{**}(\lambda,\mu)$ subject to $a \lambda + b \mu \le 1$, where constants $a>0$ and $b>0$ describe the unit rates at which self-propulsion and synchronization consume some common resource. We provide simulation results illustrating that this approach works well, when the number of particles is large. 

\subsection{Previous work and discussion of our results}

The literature on the mean-field asymptotic behavior of large-scale particle systems is vast. We will only describe the previous work that is most closely related -- a reader may find further literature reviews in the works cited here.

An MFL is a function $(f_x(t), ~x\in \R, ~t\ge 0)$, where $(f_x(t), ~x\in \R)$ is the CDF of particle locations on the real line at time $t$, in the $n\to\infty$ limit. Thus, it describes ``macro-behavior'' of a particle system, while ``micro-behavior'' often refers to the evolution of the system when the number $n$ of particles is finite. 

It is well known that the classical Kolmogorov-Petrovskii-Piscounov equation \cite{kpp-1937}
\beql{eq-kpp}
\frac{\partial}{\partial t} f_x(t) = \frac{1}{2} \frac{\partial^2 }{\partial x^2} f_x(t) +F(f_x(t)),
\eeql
describes MFLs (in our terminology) of many large-scale particle systems (cf. \cite{Man-Sch-2005, jonck-fpp-2020}
and references therein). In \eqn{eq-kpp}, $F(\cdot)$ is a non-positive function within a certain class (for example, $F(y) = -y(1-y)$ or $F(y) = -y(1-y^2)$). The seminal KKP paper \cite{kpp-1937} studies equation \eqn{eq-kpp} and proves, in particular, that equation \eqn{eq-kpp} has traveling wave solutions with any speed $v \ge v_{**}$, where the minimum speed $v_{**}$ depends on function $F(\cdot)$; it is also proved that $v_{**}$ is the average speed of the solution with initial condition being Dirac distribution concentrated at $0$ (we refer to such a solution as BMFL).
It has since been found (cf. \cite{bramson-1986, brunet-derrida-1997,brunet-derrida-2001,jonck-fpp-2020} and references therein) that many instances of particle systems, having KPP as its MFL (or, in fact, a different asymptotic limit), are such that the steady-state particle system speed $v_n$ is such that $\lim_{n\to\infty} v_n = v_{**}$. In other words, the micro-behavior is such that the system ``selects the minimum speed of the traveling waves that can arise in macro-behavior.'' This is sometimes referred to as a minimum speed selection principle (MSSP) \cite{brunet-derrida-1997,brunet-derrida-2001}.
In this paper we show that the same basic phenomenon holds for our system, even though its macro-behavior (i.e., the MFL dynamics equation), even in the special case of exponential jumps, is {\em not within KPP class}. Our proof of the traveling wave existence results for exponential jumps, while follows the general approach of KPP \cite{kpp-1937} 
(i.e., analysis of the phase portrait of a 2-dimensional ODE) require new observations -- in particular, the use of parabolic (as opposed to linear) barriers.

The first paper to formally prove an MSSP for a particle system with macro-behavior within the KKP class is \cite{bramson-1986}. Recent paper \cite{jonck-fpp-2020} considers a particle system which is just like ours, but the independent movement of each particle is a Brownian motion (instead of forward jumps); the paper proves that the MFLs of the system are within KKP class, and proves the MSSP for it. The MSSP is not limited to particle systems, whose macro-behavior within KPP class. For example, paper \cite{durrett-2011} proves MSSP for a branching selection process, whose MFL is described by an integro-differential equation with free left boundary; the paper proves, in particular, a traveling wave existence for all speeds greater or equal a critical speed. We note that the auxiliary systems with moving left and right boundaries are substantially different from the moving free boundary arising in MFL dynamics in \cite{durrett-2011}: our boundaries are ``reflecting,'' rather than ``absorbing.''

The proofs of the $\lim v_n = v_{**}$ in \cite{durrett-2011, jonck-fpp-2020} (see also references therein) rely on system monotonicity properties and on the relation with a corresponding branching process: either a branching brownian motion or a branching random walk (BRW). In fact, the critical speed $v_{**}$ is the average speed of advance of the leading (right-most) particle of the branching process. Our proof of $\limsup v_n = v_{**}$ is essentially same as in \cite{durrett-2011, jonck-fpp-2020}, except the corresponding BRW in our case is different; our proof of $\liminf v_n = v_{**}$ has relation to the proof of corresponding property in \cite{durrett-2011}, but is different. 

As we just mentioned, the critical speed $v_{**}$ is the average speed of advance of the leading (right-most) particle of the branching process. McKean \cite{mckean-1975-kpp-branching} found that BMFL for a KPP equation \cite{kpp-1937}
is such that $(f_x(t), ~x\in \R)$ is the distribution of the leading particle location for the associated branching Brownian motion with the initial particle location at $0$; we need and provide analogous BMFL characterization for our system in terms of the corresponding BRW. For a general discrete-time BRW, the average speed $v_{**}$ of the leading particle 
(as well as that of the trailing particle) was found by Biggins \cite{biggins-1976} (who extended the results of Kingman \cite{kingman-1975}); the average speed $v_{**}$ of the leading particle of the continuous-time BRW in our case follows from the results of  \cite{biggins-1976}.

The line of work  \cite{larson-1978,freidlin-1979,freidlin-1985,veret-1999} is concerned with the speed of the solutions (MFLs) of the KPP equation \cite{kpp-1937}, depending on the initial condition or, more precisely, the exponent of the right tail of the initial condition. The key results show that any speed $v \ge v_{**}$ can be achieved. In \cite{larson-1978} it is done via direct analysis of KPP solutions \cite{kpp-1937}, while \cite{freidlin-1979,freidlin-1985,veret-1999} use Feynmann-Kac representation of a solution and large deviations approach. We obtain analogous results for our system, under the additional assumption of exponential jumps, but our approach is completely different -- we prove the existence of traveling waves, including waves with moving boundaries, and use these traveling waves as lower and upper bounds of MFLs.
We note that our approach can be applied to KPP solutions \cite{kpp-1937} as well, because the existence of 
traveling waves with moving boundaries for KPP can be obtained from \cite{kpp-1937} the same way as we do for our model. In particular, we believe this gives alternative proofs of many of the MFL speed results in  \cite{larson-1978,freidlin-1979,freidlin-1985,veret-1999}, and should be applicable to other models as well.

In terms of motivation, this paper is closely related to the work on models, where a degree of ``synchronization'' is desired, meaning that the collection of particles remains ``tight,'' specifically most particles remain within $O(1)$ distance of each other uniformly in $n$. This includes the ``rank based'' model in \cite{GSS96, GMP97, St2020-wave, St2022-wave-lst},
where the instantaneous rate at which a particle jumps forward is a non-increasing function of the particle's location quantile within the empirical distribution of all particles' locations. 
The results of \cite{GSS96, GMP97, St2020-wave, St2022-wave-lst} establish most of the properties of interest: convergence to an MFL as $n\to\infty$; convergence of an MFL to a traveling wave as time goes to infinity; convergence of stationary distributions to Dirac distribution, concentrated on a traveling wave. Other models of 
distributed synchronization include those in, e.g., \cite{Man-Sch-2005, Mal-Man-2006, Manita-2009, Manita-2006, Malyshkin-2006,Manita-2014,balazs-racz-toth-2014}. 
Paper \cite{balazs-racz-toth-2014} considers a ``barycentric'' model, where the 
instantaneous rate at which a particle jumps forward is a non-increasing function of the particle's location displacement from the mean of the empirical distribution of all particles' locations; it establishes (under certain conditions) the convergence to the MFL results and finds traveling wave forms in some special cases. 
From a technical point of view, our proofs of convergence to MFL and its uniqueness closely follow the approach used in 
\cite{balazs-racz-toth-2014, St2022-wave-lst}. Paper \cite{St2020-wave} employs an artificial system with moving ``reflecting'' boundaries as a tool for analysis of the original system; this approach is one of the tools used in the present paper as well.

Finally, we note that papers \cite{brunet-derrida-1997,brunet-derrida-2001} provided first explanation of the phenomenon of ``slow convergence'' of steady-state speeds, $v_n \to v_{**}$, for particle systems having KPP solutions as MFLs. The first formal proof of this phenomenon was given, for a closely related particle system, in \cite{berard-2010}. Our simulations indicate that this phenomenon occurs in our system as well; proving this rigorously may be a subject of future work.

\subsection{Outline of the rest of the paper.}

Section~\ref{sec-notation} gives basic notation used throughout the paper.
In Section~\ref{sec-proc-def} we formally define the system process and describe its basic properties,
while Section~\ref{sec-prelim-monotone} specifically focuses on the important monotonicity properties.
Section~\ref{sec-mfl} presents our results on the convergence to -- and properties of -- mean-field limits (MFL).
Section~\ref{sec-trav-waves} defines traveling waves. In Section~\ref{sec-min-speed-principle} we 
formally define the critical speed $v_{**}$ via a minimum speed selection principle. 
In Section~\ref{sec-brw} we define the associated BRW, give McKean-type characterization of the original system's BMFL
in terms of this BRW, and show that $v_{**}$ is equal to the average speed of the leading particle of the BRW; finally, we present
%Section~\ref{sec-speed-limit} presents 
the $\lim v_n = v_{**}$ result.
In Section~\ref{sec-tws-exp} we state the traveling wave existence results in the special case of exponentially distributed jumps. In Section~\ref{sec-mfl-speed-results} we present the results on the MFL speeds: the lower bound $v_{**}$,
which holds for general jump size distribution, and speed bounds for exponential jumps case, which depend on the right tail exponent of the initial state. 
In Section~\ref{sec-conclusions} we state some further conjectures and present simulation results.
Sections~\ref{sec-MFL-B}-\ref{sec-proof-mfl-speed-exp} contain proofs.

\section{Basic notation}
\label{sec-notation}

The set of real numbers is denoted by $\R$, and is viewed as the usual Euclidean space.
As a measurable space, $\R$ is endowed with Borel $\sigma$-algebra. % $\mathcal B(\R)$. 
For scalar functions $g(x)$ of a real $x$: $\|g\|_1 = \int_x |g(x)|dx$ is $L_1$-norm; 
$g(x)$ is called $c$-Lipschitz if it is Lipschitz with constant $c \ge 0$. 
Let $\cc_b$ be the set of continuous bounded functions on $\R$, which are constant outside a closed interval (one constant value to the ``left'' of it, and possibly another constant value to the ``right'' of it.)

For functions $g(x)$ of a real $x$: 
$g(x+)$ and $g(x-)$ are the right and left limits; 
a function $g(x)$ is RCLL if it is right-continuous and has left limits at each $x$. 

A function $g$ of $x$ may be written as either $g(x)$ or $g_x$.
Notation $d_x g(x,t)$ for a multivariate function $g(x,t)$, where $x\in \R$, 
 means the differential in $x$.
 
Denote by $\cm$ the set of scalar RCLL non-decreasing
 functions $f=(f(x), ~x\in \R)$, which are (proper) probability distribution functions,
i.e., such that $f(x) \in [0,1]$, $\lim_{x\downarrow -\infty} f(x) = 0$ and $\lim_{x\uparrow \infty} f(x) = 1$.
For elements $f \in \cm$ we use the terms {\em distribution function} and {\em distribution} interchangeably.
Space $\cm$ is endowed with Levy-Prohorov metric (cf. \cite{Ethier_Kurtz}) and the corresponding topology
of weak convergence (which is equivalent to the convergence at every point of continuity of the limit);
the weak convergence in $\cm$ is denoted $\stackrel{w}{\rightarrow}$. Subspaces of $\cm$ defined later inherit its metric. 
The 
%COM LCRL version ENDCOM
inverse ($\nu$-th quantile) 
of $f \in \cm$ is $f_\nu^{-1} \doteq \inf\{y~|~f_y \ge \nu\}$, $\nu\in [0,1]$; $f_1^{-1} =\infty$ when 
$f_y < 1$ for all $y$; we will also use notation $q_\nu(f) \doteq f_\nu^{-1} $.
For a given $\nu\in (0,1)$, denote $\mathring \cm = \{f \in \cm ~|~ f_\nu^{-1} =0\}$;
the parameter $\nu$, used in the definition of $\mathring \cm$, will be clear from the context. 

Unless explicitly specified otherwise, we use the following conventions regarding random elements and random processes.
A measurable space is considered equipped with a Borel $\sigma$-algebra, induced by the 
metric which is clear from the context. A random process $Y(t), ~t\ge 0,$ always takes values in a complete separable metric space (clear
from the context), and has RCLL sample paths.
For a random process $Y(t), ~t\ge 0,$ we denote by $Y(\infty)$ the random value of $Y(t)$ in a stationary regime (which will be clear from the context). Symbol $\Rightarrow$ signifies convergence of random elements in distribution; $\stackrel{\pr}{\longrightarrow}$ means convergence in probability.
 {\em W.p.1} or {\em a.s.} means {\em with probability one.}
{\em I.i.d.} means {\em independent identically distributed.}
 For a condition/event $A$, $\bI\{A\}=1$ if $A$ holds, and $\bI\{A\}=0$ otherwise.

Space $D([0,\infty), \R)$ [resp. $D([0,\infty), \cm)$] is the Skorohod space of RCLL functions on $[0,\infty)$ taking values in 
$\R$ [resp. $\cm$], with the corresponding Skorohod ($J_1$) metric and topology (cf. \cite{Ethier_Kurtz});
$\stackrel{J_1}{\rightarrow}$ denotes the convergence in these spaces.

For a distribution $f \in \cm$ and scalar function $h(x), x \in \R$, 
$f h \doteq \int_{\R} h(x) df_x$. In some cases it is important whether or not an integration includes end points of the integration interval; 
$\int_{a+}^b h(x) df_x$ means that the integration excludes left end point $a$, 
and $\int_{a}^{b+} h(x) df_x$ means that the integration includes right end point $b$.  

For scalar functions $h(x), x \in \cx$, with some domain $\cx$, $\|h\| = \sup_{x\in \cx} |h(x)|$ is the sup-norm.
When $G_k, G$ are operators mapping the space of such functions into itself, 
$\lim G_k h = Gh$ and $G_k h \to Gh$  mean the uniform convergence: 
$\|G_k h - Gh\| \to 0$.

Suppose we have a Markov process with state space $\cx$ and transition function $P^t(x,H),$ $t\ge 0$. 
$P^t$, as an operator, is $P^t h(x)  \doteq \int_y P^t(x,dy) h(y)$, where $h$ is a scalar function with domain $\cx$;
 $I=P^0$ is the identity operator.
The process (infinitesimal) generator $B$ is 
$$
Bh \doteq \lim_{t\downarrow 0} \frac{1}{t} [P^t - I] h.
$$
Function $h$ is within the domain of the generator $B$ if $Bh$ is well-defined.

RHS and LHS mean right-hand side and left-hand side, respectively; WLOG means {\em without loss of generality}.
Abbreviation {\em w.r.t.} means {\em with respect to};
{\em a.e.} means {\em almost everywhere w.r.t. Lebesgue measure}.

We denote by $\bar J(x) = 1- J(x)$ the complementary CDF of the jump size distribution. Heaviside step function (or, Dirac distribution concentrated at point $0$) is denoted $\chi=(\chi_x = \bI\{x \ge 0\}) \in \cm$.

\section{System process definition}
\label{sec-proc-def}

We already noted that, WLOG, $\E Z= 1$; otherwise, this condition
is achieved by rescaling space. Note also that, WLOG, we can and do assume $\mu=1$; otherwise, we can achieve this condition by rescaling time.

Let $f^n(t) =(f_x^n(t), ~x\in \R)$ be the (random) empirical distribution of the 
particle locations at time $t$; namely, $f_x^n(t)$ is the fraction of particles located in $(-\infty,x]$ at time $t$. 
Clearly, $f^n(\cdot)$ is a Markov process with the state space $\cm^{(n)} \subset \cm$.
Let us also consider a centered version of the process. Fix a number $\nu \in (0,1)$,
and denote
$\mathring g_x = g_{x+q_{\nu}(g)}, ~x \in \R.$ Therefore, $\mathring g$ is the version of $g$, centered at its 
$\nu$-th quantile. Then $\mathring f^n(t)$ -- the centered version of $f^n(t)$ -- is also a Markov process,
with state space $\mathring \cm^{(n)} = \cm^{(n)} \cap \mathring \cm$. 
This process is regenerative, with the regenerations occurring at time points when all particles ``collapse'' to one location.
A regeneration cycle has finite mean; indeed, for any time interval of a fixed length $\epsilon>0$,
with probability at least some $\delta=\delta(\epsilon,n)>0$,
 at the end of the interval all particles will be at a single location. 
Therefore the process is positive recurrent,
with unique stationary distribution. Denote by $\mathring f^n(\infty)$ a random element with the distribution equal to the stationary distribution of $\mathring f^n(\cdot)$.

Consider a version $f^n(\cdot)$ of the system process, such that the centered process $\mathring f^n(\cdot)$
is stationary. Specifically, consider the process $f^n(\cdot)$ with $f^n(0)$ equal in distribution to $\mathring f^n(\infty)$.
Then, for each $n$, the steady-state average speed $v_n$ of the particle system is finite and given by
$$
v_n = \frac{1}{T} \E [q_\beta(f^n(t)) - q_\beta(f^n(0))],
$$
where $\beta\in (0,1)$ and $t>0$ can be chosen arbitrarily.
Also, $v_n$ is the steady-state
average speed of the center of mass, i.e. of the distribution $f^n(t)$ mean: 
$$
v_n = \lambda + \E \int_{-\infty}^{\infty} d_y \mathring f_y^n(\infty) \int_y^\infty  d_\xi \mathring f_{\xi}^n(\infty) (\xi-y).
$$
Given the regenerative structure of the (centered) process, the steady-state speed $v_n$ can also be defined 
in terms of probability one convergence, regardless of the initial state. In particular, 
for any initial state of the system,
$$
v_n = \lim_{t\to\infty} \frac{1}{t} [D^n(t) - D^n(0)], ~~w.p.1,
$$
where $D^n(t)$ is the location of the leading (right-most) particle. 

We remark that, trivially, $\lambda=0$ implies $v_n=0$, because then w.p.1 all particles will end up at a single location and will stop moving after that. Therefore, $v_n > 0$ if and only if $\lambda>0$, and we have a trivial lower bound $v_n \ge \lambda$ for any $n$.

Recall that we assumed $\mu=1$. If we do need expressions for a case when $\mu\ne 1$, they are obtained 
form the corresponding expressions for the $\mu=1$ case. For example, 
if $v_n=v_n(\lambda,\mu)$ is the steady-state speed as a function of $\lambda$ and $\mu>0$, then
$v_n(\lambda,\mu) = \mu v_n(\lambda/\mu,1)$.

\section{Monotonicity properties}
\label{sec-prelim-monotone}

Order relation (stochastic dominance) $g^{(1)} \preceq g^{(2)}$ between two distributions $g^{(1)}, g^{(2)} \in \cm$  means $g^{(1)}_x \ge g^{(2)}_x, ~\forall x$.

Process $f^n(\cdot)$ possesses some simple monotonicity properties. 
If we have two versions of this process $f^n(\cdot)$, labeled $f^{n,(1)}(\cdot)$ and $f^{n,(2)}(\cdot)$,
with fixed initial conditions
$f^{n,(1)}(0) \preceq f^{n,(2)}(0)$, then these two processes can be coupled (constructed on a common probability space) so that $f^{n,(1)}(t) \preceq f^{n,(2)}(t)$
at all times $t\ge 0$. Indeed, it suffices to set a one-to-one correspondence between particles in the two systems,
such that each particle in the second system in initially ahead of its corresponding particle in the first system,
and consider a probability space such that each pair of the corresponding particles has same realizations
of the time instants and sizes of the independent jumps and same realizations of the time instants and ``targets'' of the synchronization  jumps.
Then, clearly, each particle in the second system will remain ahead of its counterpart in the first system at all times.

This property easily generalizes to the case when the second (``larger'') process has larger number of particles, i.e.
we compare $f^{n,(1)}(\cdot)$ and $f^{k,(2)}(\cdot)$, where $n \le k$. In this case, we will use the order relation
$f^{n,(1)}(t) \preceq_l f^{k,(2)}(t)$, which denotes the situation when $f^{n,(1)}(t) \preceq f^{n(k),(2)}(t)$,
where $f^{n(k),(2)}(t)$ is the empirical distribution of the leading $n$ particles of the state $f^{k,(2)}(\cdot)$.
Again, if $f^{n,(1)}(0) \preceq_l f^{k,(2)}(0)$, then the processes can be coupled so that $f^{n,(1)}(t) \preceq_l f^{k,(2)}(t)$ prevails at all times. This extension has been used in the literature in different contexts (cf the proof of theorem 2.3(4) in \cite{jonck-fpp-2020}, and references therein), but is also easy to see directly. Indeed, the synchronization jumps can be equivalently defined as follows: a particle gets a ``synchronization urge'' as a Poisson process of rate $\mu$ (rescaled to $1$),
and it relocates another particle, chosen uniformly at random, to its own location if the chosen particle happens to be behind it. Given this interpretation, the coupling construction is straightforward.

Since the process is such that (a) for any initial state, $v_n = \lim_{t\to\infty} D^n(t)/t$ a.s., where $D^n(t)$ is the location of leading particle at time $t$, and (b) monotonicity w.r.t. $\preceq_l$ holds, we see that the steady-state speed $v_n$ is non-decreasing in $n$:
\beql{eq-vn-increasing}
v_n \le v_{n+1}, ~~\forall n \ge 1.
\eeql
(Again, this property was used in the proof of theorem 2.3(4) in \cite{jonck-fpp-2020} for a different system, but it obviously holds for any system satisfying (a) and (b).)

The monotonicity with respect to $\preceq$ (and, more generally, $\preceq_l$) 
further easily generalizes in several directions. For example, it still holds if the ``larger'' process 
$f^{k,(2)}(\cdot)$ has larger parameter $\lambda$; and/or if it is ``helped'' by a left boundary (can be static or moving in any way); and/or if the ``smaller'' process $f^{n,(1)}(\cdot)$  is ``impeded'' by a right boundary (can be static of moving
in any way); etc. We will not give formal statements for all these generalizations/variations, and will just refer to monotonicity properties described in this section as {\em monotonicity.}

\section{Mean-field limits}
\label{sec-mfl}

In this section we state our results on the convergence to and characterization of the mean-field limits.
The results are for the original system, the system with fixed right boundary (which is primarily a tool for the analysis
of the original system), and for the two systems with moving right and left boundaries, respectively.
We will use the same notations $f^n(\cdot)$ and $f(\cdot)$ for the process and a mean-field limit, respectively, for each of the four systems.
It should be clear that they do depend on the system, and when we use them later in the paper it will be clear from the context which system they refer to.

\subsection{Original system}

The following definition of a mean-field model describes what is natural to expect a mean-field limit to be. 
Theorem~\ref{th-finite-interval} then shows that this is indeed the case.

\begin{definition}
\label{def-mfm}
A function $f_x(t), ~x\in \R, ~t\in \R_+,$ will be called a {\em mean-field model (MFM)} 
if it satisfies the following conditions.\\
 (a) For any $t$, $f(t) = (f_x(t), x\in \R) \in \cm$.\\
(b) For any $x$, $f_x(t)$ is non-increasing $(\lambda+1)$-Lipschitz in $t$.\\
(c) For any $x$, for any $t$ where the partial derivative
$(\partial/\partial t) f_x(t)$ exists (which is almost all $t$ w.r.t. Lebesgue measure, by the Lipschitz property), equation 
\beql{eq-dyn}
\frac{\partial}{\partial t} f_x(t) = -\lambda \int_{-\infty}^x d_y f_y(t) (1-J(x-y)) - f_x(t) (1-f_x(t)).
\end{equation}
holds. 
\end{definition}

Denote by $L^{(n)}$ the generator of the process
$f^n(\cdot)$. 
For any $h\in\cc_b$, function $f^n h$ of element $f^n \in \cm^{(n)}$ is within the domain of $L^{(n)}$ (where we use the fact each function in $\cc_b$ is constant outside a closed interval), and
$$
L^{(n)} [f^n h] = 
\lambda \E_Z \int_{-\infty}^{\infty} d_y f_y^n (h(y+Z)-h(y))
+ \int_{-\infty}^{\infty} d_y  f_y^n \int_y^\infty  d_\xi  f_{\xi}^n (h(\xi)-h(y)),
$$
where the expectation is over the distribution of the random jump size $Z$.  
We also formally define the ``limit'' $L$ of $L^{(n)}$ (for elements $f \in \cm$) as
$$
L [f h] = 
\lambda \E_Z \int_{-\infty}^{\infty} d_y f_y (h(y+Z)-h(y))
+ \int_{-\infty}^{\infty} d_y  f_y \int_y^\infty  d_\xi  f_{\xi} (h(\xi)-h(y)).
$$

\begin{thm}
\label{th-finite-interval}
Suppose $f^n(0) \stackrel{w}{\rightarrow} f(0)$, where $\{f^n(0)\}$ is deterministic sequence of $f^n(0) \in \cm^{(n)}$,
and $f(0) \in \cm$. 
Then $f^n(\cdot) \Rightarrow  f(\cdot)$ in $D([0,\infty), \cm)$, where $f(\cdot) \in D([0,\infty), \cm)$ is deterministic, uniquely determined by $f(0)$. Moreover, $f(\cdot)$ is a continuous element of $D([0,\infty), \cm)$, which satisfies 
\beql{eq-pde-oper}
f(t) h - f(0) h - \int_0^t L f(s) h ds =0, ~~\forall h\in \cc_b, ~\forall t\ge 0.
\eeql
This trajectory $f(\cdot)$ is called the {\em mean-field limit (MFL)} with initial state $f(0)$.
Furthermore, MFL $f(\cdot)$ is the unique mean-field model (see Definition~\ref{def-mfm}) with initial state $f(0)$.
\end{thm}

The proof of Theorem~\ref{th-finite-interval} is in Section~\ref{sec-MFL-orig}.

The MFL $f(\cdot)$ with initial state $f(0)=\chi$ (i.e., the Dirac distribution concentrated at a single point $0$) we will call the {\em benchmark MFL (BMFL)}.

\subsection{System with fixed right boundary}

As an intermediate step towards the proof of Theorem~\ref{th-finite-interval}, we will prove its special case,
Theorem~\ref{th-finite-interval-B} below, which applies to the system with fixed right boundary, defined as follows.

Note that, for any fixed $B$, the evolution of $f_x^n(\cdot), x \in (-\infty, B)$ in our original system is independent of the actual locations and/or evolution of particles in $[B, \infty)$. Therefore, if we are only interested in the evolution of $f_x^n(\cdot), x \in (-\infty, B)$,  we may as well assume that any particle in $[B, \infty)$ is located exactly at $B$ (that is, 
$f_B^n(t)=1$ at all times). In other words, we can assume that there is a fixed right boundary point $B$, and
any particle which ``tries to jump over $B$'' lands (absorbed) at exactly $B$. This is the modified system 
(with fixed right boundary) that we consider here. We see that the process $f^n(\cdot)$ for this system is a  projection 
of the original process -- it is such that $f^n(t) \in \cm^B \doteq \{g \in \cm ~|~ g_B=1\}$ for all $t$.

Denote by $L^{(n),B}$ the generator of the process
$f^n(\cdot)$ for this system. 
For any $h\in\cc_b$, function $f^n h$ of element $f^n \in \cm^{(n)} \cap \cm^B$ is within the domain of $L^{(n),B}$,
and
$$
L^{(n),B} [f^n h] = 
\lambda \E \int_{-\infty}^{B} d_y f_y^n (h((y+Z) \vee B)-h(y))
+ \int_{-\infty}^{B} d_y  f_y^n \int_y^{B+}  d_\xi  f_{\xi}^n (h(\xi)-h(y)).
$$
where the expectation is over the distribution of the random jump size $Z$.  
We also formally define the ``limit'' $L^{B}$ of $L^{(n),B}$ (for elements $f \in \cm \cap \cm^B$) as
$$
L^{B} [f h] = 
\lambda \E \int_{-\infty}^{B} d_y f_y (h((y+Z) \vee B)-h(y))
+ \int_{-\infty}^{B} d_y  f_y \int_y^{B+}  d_\xi  f_{\xi} (h(\xi)-h(y)).
$$

Observe the following.  For a fixed $h\in\cc_b$, if $B$ is large enough so that $h(x)$ is constant for 
all $x \ge B$, then $L^{(n),B} [f^n h]= L^{(n)} [f^n h]$ and $L^{B} [f^n h]= L [f^n h]$, where 
$L^{(n)}$ and $L$ are the operators we defined for the original system.

\begin{thm}
\label{th-finite-interval-B}
Suppose $f^n(0) \stackrel{w}{\rightarrow} f(0)$, where $\{f^n(0)\}$ is deterministic sequence of $f^n(0) \in \cm^{(n)} \cap \cm^B$, 
and $f(0) \in \cm^B$.  
Then $f^n(\cdot) \Rightarrow  f(\cdot)$ in $D([0,\infty), \cm^B)$, where $f(\cdot) \in D([0,\infty), \cm^B)$ is deterministic, uniquely determined by $f(0)$.  
Moreover, $f(\cdot)$ is a continuous element of $D([0,\infty), \cm^B)$, which satisfies 
\beql{eq-pde-oper-B}
f(t) h - f(0) h - \int_0^t L^B f(s) h ds =0, ~~\forall h\in \cc_b, ~\forall t\ge 0.
\eeql
This trajectory $f(\cdot)$ is called the {\em mean-field limit (MFL)} with initial state $f(0)$.
Furthermore, MFL $f(\cdot)$ is the unique mean-field model with initial state $f(0)$; 
the mean-field model here is defined as in Definition~\ref{def-mfm} except \eqn{eq-dyn} holds for $x < B$,
and $f_x(t) \equiv 1$ for $x\ge B$ and all $t$.
\end{thm}

The proof of Theorem~\ref{th-finite-interval-B} is in Section~\ref{sec-MFL-B}.

\subsection{System with moving right boundary}

Consider now the system with the right boundary point $B$ moving right at constant speed $v>0$, i.e. $B=B_0+vt$
for some constant $B_0$. Just like in the system with fixed right boundary, any particle which ``tries to jump over the boundary'' lands exactly at its current location. However, particles landing at the right boundary are not ``absorbed'' in it,
as the boundary keeps moving right while the particles do not move until they jump again. If $B$ is the boundary location
at time $t$, then the state $f^n(t) \in \cm^B$. Formally speaking, the state of the Markov process has the additional component -- real number $B$; somewhat abusing notation, we will not include this component into the state descriptor.

It is important to note that the process with moving right boundary is {\em not} a projection of the original process.

Denote by $L^{(n), \{r\}}$ the generator of the process
$f^n(\cdot)$ for this system. 
For any $h\in\cc_b$, function $f^n h$ of element $f^n$ (or rather pair $(f^n,B)$), satisfying
$f^n \in \cm^{(n)} \cap \cm^B$, is within the domain of $L^{(n), \{r\}}$, and
$$
L^{(n), \{r\}} [f^n h] = 
\lambda \E \int_{-\infty}^{B} d_y f_y^n (h((y+Z) \vee B)-h(y))
+ \int_{-\infty}^{B} d_y  f_y^n \int_y^{B+}  d_\xi  f_{\xi}^n (h(\xi)-h(y)),
$$
where the expectation is over the distribution of the random jump size $Z$.  
We also formally define the ``limit'' $L^{\{r\}}$ of $L^{(n), \{r\}}$ (for pairs $(f,B)$ satisfying
$f \in \cm^B$) as
$$
L^{\{r\}} [f h] = 
\lambda \E \int_{-\infty}^{B} d_y f_y (h((y+Z) \vee B)-h(y))
+ \int_{-\infty}^{B} d_y  f_y \int_y^{B+}  d_\xi  f_{\xi} (h(\xi)-h(y)).
$$
Operators $L^{(n), \{r\}}$ and $L^{\{r\}}$ are same as operators $L^{(n), B}$ and $L^B$, respectively, for the system with fixed right boundary at $B$, but here $B$ is part of the process state, rather than a fixed constant.

\begin{thm}
\label{th-finite-interval-drift-1}
Let $v>0$ and $B_0$ be fixed, and recall that boundary $B=B_0 +vt$.
Suppose $f^n(0) \stackrel{w}{\rightarrow} f(0)$, where $\{f^n(0)\}$ is deterministic sequence of $f^n(0) \in \cm^{(n)} \cap \cm^{B_0}$, 
and $f(0) \in \cm^{B_0}$.  
Then $f^n(\cdot) \Rightarrow  f(\cdot)$ in $D([0,\infty), \cm)$, where $f(\cdot) \in D([0,\infty), \cm)$ is deterministic, uniquely determined by $f(0)$; $f^n(t) \in \cm^{(n)} \cap \cm^{B}$ and $f(t) \in \cm^{B}$ for all $t$.
Moreover, $f(\cdot)$ is a continuous element of $D([0,\infty), \cm)$, which satisfies 
\beql{eq-pde-oper-drift-1}
f(t) h - f(0) h - \int_0^t L^{\{r\}} f(s) h ds =0, ~~\forall h\in \cc_b, ~\forall t\ge 0.
\eeql
This trajectory $f(\cdot)$ is called the {\em mean-field limit (MFL)} with initial state $f(0)$.
Furthermore, MFL $f(\cdot)$ is the unique mean-field model with initial state $f(0)$; 
the mean-field model here is defined as in Definition~\ref{def-mfm} except \eqn{eq-dyn} holds for $x < B (=B_0 + vt)$.
\end{thm}

The proof of Theorem~\ref{th-finite-interval-drift-1} is in Section~\ref{sec-MFL-1}.

\subsection{System with moving left boundary}

Consider yet another modification of our original system, namely the system with left boundary 
$A$ moving right at constant speed $v>0$, i.e. $A=A_0+vt$ for some constant $A_0$. 
In this system there are no particles to the left of $A$ -- as the boundary moves right it ``drags forward'' any particle 
that it encounters and such particle stays at the (moving) boundary until it jumps forward. 
The corresponding process is {\em not} a projection of our original process. 

By the process definition, there are no particles located to the left of moving boundary $A$. 
However, for the purposes of analysis it will be convenient to adopt an equivalent view of the process,
which allows particles to be to the left of $A$. Specifically, assume that particles make synchronization jumps as 
in the original system, but when a particle located at $y$ makes an independent jump, it first instantly moves to point $y \vee A$ and then jumps.
Clearly, this process is such that the evolution in time of $(f^n_x(t), x \ge A)$ is exactly same as it would be
if the moving left boundary would drag the particles that it encounters with it. So, this is how we will 
define the process $f^n(\cdot)$ for this system. Then $f^n(t) \in \cm^{(n)}$ for all $t$, and boundary location $A$ is implicitly a part of the process definition.

Denote by $L^{(n), \{l\}}$ the generator of the process
$f^n(\cdot)$. 
For any $h\in\cc_b$, function $f^n h$ of element $f^n$ (or rather pair $(f^n,A)$), satisfying
$f^n \in \cm^{(n)}$ is within the domain of $L^{(n), \{r\}}$, and
$$
L^{(n), \{l\}} [f^n h] = 
\lambda \E \int_{-\infty}^{\infty} d_y f_y^n ( h((y \vee A)+Z)-h(y) )
+ \int_{-\infty}^{\infty} d_y  f_y^n \int_y^{\infty}  d_\xi  f_{\xi}^n (h(\xi)-h(y)).
$$
where the expectation is over the distribution of the random jump size $Z$.  
We also formally define the ``limit'' $L^{\{l\}}$ of $L^{(n), \{l\}}$  (for pairs $(f,A)$ with
$f \in \cm$) as
$$
L^{\{l\}} [f^n h] = 
\lambda \E \int_{-\infty}^{\infty} d_y f_y ( h((y \vee 0)+Z)-h(y) )
+ \int_{-\infty}^{\infty} d_y  f_y \int_y^{\infty}  d_\xi  f_{\xi} (h(\xi)-h(y)).
$$

\begin{thm}
\label{th-finite-interval-drift-2}
Let $v>0$ and $A_0$ be fixed, and recall that boundary $A=A_0 +vt$.
Suppose $f^n(0) \stackrel{w}{\rightarrow} f(0)$, where $\{f^n(0)\}$ is deterministic sequence of $f^n(0) \in \cm^{(n)}$,
and $f(0) \in \cm$.  
Then $f^n(\cdot) \Rightarrow  f(\cdot)$ in $D([0,\infty), \cm)$, where $f(\cdot) \in D([0,\infty), \cm)$ is deterministic, uniquely determined by $f(0)$. 
Moreover, $f(\cdot)$ is a continuous element of $D([0,\infty), \cm)$, which satisfies 
\beql{eq-pde-oper-drift-2}
f(t) h - f(0) h - \int_0^t L^{\{l\}} f(s) h ds =0, ~~\forall h\in \cc_b, ~\forall t\ge 0.
\eeql
This trajectory $f(\cdot)$ is called the {\em mean-field limit (MFL)} with initial state $f(0)$.
Furthermore, MFL $f(\cdot)$ is the unique mean-field model with initial state $f(0)$; 
the mean-field model here is defined as in Definition~\ref{def-mfm} except \eqn{eq-dyn}
is replaced by
\beql{eq-dyn-l2}
\frac{\partial}{\partial t} f_x(t) = -\lambda f_x(t) - f_x(t) (1-f_x(t)), ~x \le A,
\end{equation}
\beql{eq-dyn-l3}
\frac{\partial}{\partial t} f_x(t) = -\lambda  f_A(t) (1-J(x-A)) -\lambda \int_{A}^x d_y f_y(t) (1-J(x-y)) - f_x(t) (1-f_x(t)), ~x > A.
\end{equation}
\end{thm}

The proof of Theorem~\ref{th-finite-interval-drift-2} is in Section~\ref{sec-MFL-2}.

\subsection{Monotonicity of mean-field limits.}
\label{sec-mfl-monotone}

The process monotonicity properties described in Section~\ref{sec-prelim-monotone},
as well as MFL uniqueness for each initial state, imply the following monotonicity property of MFLs.

\begin{lem}
\label{lem-mfl-monotone}
For two MFLs $f^{(1)}(\cdot)$ and $f^{(2)}(\cdot)$,  $f^{(1)}(0) \preceq f^{(2)}(0)$ implies 
$f^{(1)}(t) \preceq f^{(2)}(t)$ for all $t\ge 0$.
\end{lem}

Just as the process monotonicity still holds under various generalizations (see Section~\ref{sec-prelim-monotone}),
Lemma~\ref{lem-mfl-monotone} holds under those generalizations as well.

\section{MFLs that are traveling waves}
\label{sec-trav-waves}

An MFL (MFM) 
$f(\cdot)$ for our original system is a traveling wave, if for some speed $v> 0$ and some 
$\phi = (\phi_x, x \in \R) \in \cm$, $f_x(t) =\phi_{x-vt}$, in which case we call $\phi$ a {\em traveling wave shape (TWS).} Substituting $f_x(t) =\phi_{x-vt}$ into
\eqn{eq-dyn}, we can easily see that  a TWS, {\em if it exists}, must satisfy
\beql{eq-wave}
v \phi'_x = \lambda \int_{-\infty}^x \phi'_y (1-J(x-y)) dy + \phi_x (1-\phi_x)
\end{equation}
for each $x$, and in fact the derivative $\phi'_x$ must be continuous in $x$.

Since a traveling wave is an MFL (which is a limit of particle systems dynamics), we see that {\em if a TWS with speed $v$ exists then necessarily $v \ge \lambda$.}

 {\em In the special case when $\lambda=0$}, the integro-differential equation \eqn{eq-wave} becomes 
logistic differential equation \beql{eq-wave-special}
v \psi'_x =  \psi_x (1-\psi_x).
\end{equation}
In this case, a unique (up to a shift) TWS does exist for every speed $v>0$. This follows, for example, from \cite[proposition 5.1]{Man-Sch-2005}, but is, of course, well-known:  
the relevant solution to \eqn{eq-wave-special}  
is 
\beql{eq-wave-zero1}
\psi_x = \frac{1}{1+ e^{-(1/v)(x+c)}}, ~x\in \R, ~~\mbox{where $v>0$ and $c\in \R$ are parameters}.
\end{equation}
(There is no contradiction here with the fact that for any {\em finite} $n$, $\lambda=0$ implies that the steady-state speed $v_n=0$. 
Traveling waves %, on the other hand, 
are MFLs, 
which means that we first take the $n\to\infty$ limit and then look at the time-evolution of the limit.
So, if time interval $[0,T]$ is fixed, $n$ is very large, and the initial state $f^n(0)$ is close to $\psi$ in \eqn{eq-wave-zero1}, then
the evolution of $f_x^n(t)$ in $[0,T]$ is close to the traveling wave $\psi_{x-vt}$, moving at speed $v$. This does not contradict the fact that, eventually, all particles will assemble at the location of the leading particle and will stop moving, i.e. $v_n=0$.)

{\em In the special case of exponential jump size distribution,} $J(x) =1- e^{-x}$, for any $\lambda\ge 0$, 
the integro-differential equation \eqn{eq-wave} becomes an ODE
\beql{eq-tws-exp-1}
v \phi''_x =  (1+\lambda-v-2\phi_x) \phi'_x + \phi_x (1-\phi_x).
\eeql
Indeed, in this case \eqn{eq-wave} can be written as 
$$
v \phi'_x e^x = \lambda \int_{-\infty}^x \phi'_y e^y dy + \phi_x (1-\phi_x) e^x,
$$
which, after differentiating in $x$ gives \eqn{eq-tws-exp-1}.

We will also consider MFLs %(MFMs) 
that are traveling waves for the modified systems with moving right and left boundaries; $v>0$ is the speed of the boundary. 
Clearly, in a system with a boundary moving at speed $v$, any traveling wave speed is also $v$. 
For both systems, for a given speed $v>0$ of the boundary, 
an MFL %(MFM) 
$f(\cdot)$  is a traveling wave (of speed $v$), if for some 
$\phi = (\phi_x, x \in \R) \in \cm$, $f_x(t) =\phi_{x-vt}$, where $\phi$ is called a traveling wave shape (TWS) for the corresponding system; WLOG we can assume that the initial location of the boundary is $0$. 

Substituting $f_x(t) =\phi_{x-vt}$ into
\eqn{eq-dyn}, we see that  a TWS $\phi$ for the system with moving right boundary, {\em if it exists},  must be continuous with $\phi_0=1$ and must satisfy \eqn{eq-wave}
for each $x<0$, and in fact the derivative $\phi'_x$ must be continuous in $(-\infty,0]$ (if at $x=0$ we consider left derivative). In the special case of $J(x) =1-e^{-x}$, \eqn{eq-wave} reduces to \eqn{eq-tws-exp-1}.

Substituting $f_x(t) =\phi_{x-vt}$ into \eqn{eq-dyn-l3}, we obtain that 
a TWS $\phi$ for the system with moving left boundary, {\em if it exists}, must be continuous for $x\in [0,\infty)$ and
must satisfy 
\beql{eq-wave-1}
v \phi'_x = \lambda \phi_0 (1-J(x))  + \lambda \int_{0+}^x \phi'_y (1-J(x-y)) dy + \phi_x (1-\phi_x)
\end{equation}
for each $x>0$, and in fact the derivative $\phi'_x$ must be continuous in $[0,\infty)$ (if at $x=0$ we consider right derivative). Note that the initial condition $v \phi'_0 = \lambda \phi_0 + \phi_0(1-\phi)$ must hold.
In the special case of $J(x) =1-e^{-x}$, \eqn{eq-wave-1} reduces to \eqn{eq-tws-exp-1}.

\section{Definition of the critical speed $v_{**}$ via a minimum speed selection principle.}
\label{sec-min-speed-principle}

In this section we formally define the critical value $v_{**}$ of the speed, which plays central role in our results and analysis. Assume $\lambda>0$. 
Denote by $L(s)=  \int_0^\infty e^{-sx} dJ(x)$ is the Laplace transform of the jump distribution $J(\cdot)$.
As a function of real $s$ it is, of course, a convex non-negative non-increasing function.
Denote by $\alpha$ the the tail exponent of $J(\cdot)$:
\beql{eq-alpha-def}
\alpha \doteq \liminf_{x\to\infty} \frac{-\log (1-J(x))}{x} = \sup \{\zeta \ge 0 ~|~ L(-\zeta) < \infty\} \in [0,\infty].
\eeql
Suppose, $\alpha>0$ (with the case $\alpha=\infty$ allowed). 
Then $L(s) <+\infty$ for $s> - \alpha$, $L(s) = +\infty$ for $s < - \alpha$, and 
$L(s) \uparrow L(-\alpha)$ as $s \downarrow -\alpha$; $L(-\alpha)$ may be finite or infinite, but it is necessarily infinite when $\alpha=\infty$.

The following heuristic argument, leading to \eqn{eq-v-gamma} 
is analogous to the one in, e.g., \cite{brunet-derrida-1997}, where it is applied to the KPP equation \eqn{eq-kpp}. 
This heuristic argument serves as a motivation 
for the {\em rigorous} definition of speed value $v_{**}$. Assume that a traveling wave shape $\phi = (\phi_x, ~x\in \R)$  with speed $v>0$ exists, and that {\em the front tail of $\phi$  decays exponentially:} $1-\phi_x \sim e^{-\zeta x}$, as $x\to\infty$, where $\zeta>0$. 
 Then, when $x$ is large, equation \eqn{eq-wave} ``becomes:''
\beql{eq-tail}
v \zeta e^{-\zeta x} = \lambda \int_{-\infty}^x \zeta e^{-\zeta y} \bar J(x-y) dy + e^{-\zeta x}.
\eeql
In \eqn{eq-tail} 
$$
\int_{-\infty}^x e^{-\zeta y} \bar J(x-y) dy = e^{-\zeta x} \int_0^\infty e^{\zeta \xi} \bar J(\xi) d\xi
= \frac{1}{\zeta} e^{-\zeta x} [L(-\zeta)-1],
$$
so we obtain the following dependence of speed $v$ on $\zeta>0$:
\beql{eq-v-gamma}
v(\zeta) = \frac{1}{\zeta} [\lambda L(-\zeta) - \lambda +1].
\eeql
Relation \eqn{eq-v-gamma}, viewed formally, is a starting point of the following rigorous definitions and observations.
The dependence $v(\zeta)$ is convex, because %as we saw just above, 
$(1/\zeta) [L(-\zeta)-1]$
%$$
% \frac{1}{\zeta}  [L(-\zeta)-1] 
%$$
is the Laplace transform of function $\bar J(\cdot)$ at point $s=-\zeta$. 
Note also that $v(\zeta)\to +\infty$ as $\zeta\downarrow 0$. It is easy to see that the minimum value $v_{**}$ of 
$v(\zeta)$
is attained at the unique, positive, finite value
\beql{eq-gamma-opt}
\zeta_{**} \doteq \argmin_{\zeta>0} v(\zeta) \in (0,\alpha], ~~\mbox{so that}~~ v_{**} \doteq v(\zeta_{**}) = \min_{\zeta>0} v(\zeta).
\eeql
Because $\zeta_{**} \in (0,\alpha]$ and $L(-\zeta) > 1$ for $\zeta>0$, note from \eqn{eq-v-gamma} that 
\beql{eq-speed-ge-inv-alpha}
v_{**} > 1/\alpha.
\eeql  

If $\zeta_{**} < \alpha$, then it necessarily solves equation 
\beql{eq-gamma-opt-typical}
\lambda \zeta_{**} L'(-\zeta_{**}) + \lambda L(-\zeta_{**}) -\lambda +1 = 0;
\eeql
and if $L(-\alpha) =\infty$, then necessarily $\zeta_{**} < \alpha$.
Note for future reference that the inverse function to $v(\zeta), ~0<\zeta \le \zeta_{**},$ is a continuous strictly decreasing convex function $\zeta(v), v_{**} \le v < \infty,$ with $\zeta(v) \downarrow 0$ as $v\uparrow\infty$.

In the special case $\bar J(x) = e^{-x}$, we have $\alpha=1$, 
\eqn{eq-v-gamma} becomes
\beql{eq-v-gamma-exp}
v(\zeta) = \frac{\lambda}{1-\zeta} + \frac{1}{\zeta},
\eeql
and then
$$
\zeta_{**} = \zeta_* \doteq \frac{1}{1+\sqrt{\lambda}}, ~~ 
v_{**} = v_* \doteq (1+\sqrt{\lambda})^2 .
$$
The inverse function $\zeta(v), v_{*} \le v < \infty,$ to function 
$v(\zeta), ~0<\zeta \le \zeta_{*},$ in this case is
\beql{eq-v-gamma-exp-inverse}
\zeta(v) = \frac{(1+v-\lambda) - \sqrt{(1+v-\lambda)^2 - 4v}}{2v}.
\eeql
Recall that the above definition of $v_{**}$ assumes $\alpha>0$.
In view of \eqn{eq-speed-ge-inv-alpha}, we adopt a convention that $v_{**} = +\infty$ when $\alpha=0$.
All results in this paper, involving the critical speed $v_{**}$, are valid for the case when $\alpha=0$ and $v_{**} = +\infty$.
They easily follow from the corresponding results for the case when $\alpha > 0$ and $v_{**} < \infty$,
by using monotonicity (Section~\ref{sec-prelim-monotone}). Namely, if $\alpha=0$, we can compare 
the actual process to an auxiliary one with the jump size distribution slightly changed so that it is stochastically smaller, 
with an arbitrarily small positive tail exponent $\alpha$.
The auxiliary process, which is stochastically dominated by the actual one, can have arbitrarily large value of $v_{**}$.
Thus, the corresponding results for the actual process hold with $v_{**}=\infty$.

\section{Limit of steady-state speeds and the associated branching random walk}
\label{sec-brw}

\subsection{McKean-type characterization of the original system's BMFL.}

Theorem~\ref{prop-mckean-type}, presented in this section, and its proof are analogous to those in \cite{mckean-1975-kpp-branching}, for the system where the independent movement of each particle is a Brownian motion.

Denote by $W(t), ~t \ge 0,$ a single particle independent jump (Markov) process, starting at some point $y$ at initial time $t=0$. Specifically, $W(t)$ is the location at time $t$ of the particle that starts at $0$, and makes i.i.d. jumps,
distributed as $J(\cdot)$, at time points of a Poisson process of rate $\lambda>0$. 

Consider the following process, which we will call {\em the associated (with our particle system) branching random walk (BRW).}
The process starts with a single particle, located initially at $0$. Its 
location $W_1(t)$ evolves a single particle independent jump process. Any particle present in the system generates 
a new particle in the same location (i.e. splits into two) according to a unit rate Poisson process. Each newly created particle continues as a single particle independent jump process. We label particles in the order of their creation, so that
at time $t$ their locations are $W_1(t), \ldots, W_{N(t)}(t)$, where $N(t)$ is the total number of particles. 
Let $D(t) \doteq \max_i W_i(t)$ be the location of the right-most particle of the associated BRW.

\begin{thm}
\label{prop-mckean-type}
The associated BRW is such that
$$
\pr\{D(t) \le x\} = f_x(t),
$$
where $f(t)$ is the BMFL.
\end{thm}

The proof of Theorem~\ref{prop-mckean-type} is in Section~\ref{sec-mckean}.

\subsection{Value $v_{**}$ is the limiting average speed of the leading particle of the associated BRW.}
\label{sec-brw-speed}

The following proposition can be obtained from a general result on the average speed of progress of a leading particle of a branching random walk in {\em discrete-time}, namely theorem 3 in \cite{biggins-1976}. 
\begin{prop}
\label{prop-speed-brw}
Consider the associated BRW, starting with one particle. Let $D(t)$ be the location of the leading particle at time $t$. Then, 
\beql{eq-speed-brw}
D(t)/t \to v_{**}, ~~~\mbox{w.p.1}.
\eeql
\end{prop}

Indeed, the process can viewed as follows. Each particle already present in the process ``waits'' for an independent, exponentially distributed time with mean $1/(1+\lambda)$, and then either jumps forward with probability $\lambda/(1+\lambda)$ or splits into two particles at the same location with probability $1/(1+\lambda)$. Suppose we discretize time, with the time step being small $\delta>0$.
We obtain a process ``slower'' than original one, if we assume that the time $\tau$ a particle ``waits'' until next event is geometric, with $\pr\{\tau = j\delta\} = \epsilon (1-\epsilon)^{j-1}$, $j=1,2,\ldots$, $\epsilon=1-\exp(-(1+\lambda)\delta)$.
Particles present at time $n\delta$ represent the $n$-generation of the particles.  
The process is such that, after each time step, a particle has either one descendant in the same locations with probability $1-\epsilon$, or two descendants in the same location with probability $\epsilon/(1+\lambda)$,
or one descendant located at the current particle location  + $Z$. To this discrete-time process theorem 3 in \cite{biggins-1976} applies to obtain the leading particle average speed of progress, which is a lower bound of the speed of progress for the original process. To obtain a discrete-time process, which is ``faster'' than the original one, we assume that the time $\tau$ a particle ``waits'' until next event is geometric, with $\pr\{\tau = j\delta\} = \epsilon (1-\epsilon)^{j}$, $j=0,1,2,\ldots$.
(A subtlety here is that $\tau=0$ with non-zero probability, so that the distribution of the set of one time step descendants,
i.e. after time $\delta$, has to account for that.) Applying theorem 3 in \cite{biggins-1976} to this discrete-time process, we obtain an upper bound on the leading particle average speed of progress for the original process. Letting 
$\delta\downarrow 0$, we can check that both bounds converge to exactly $v_{**}$.

\subsection{Limit of steady-state speeds}
\label{sec-speed-limit}

\begin{thm}
\label{thm-speed}
For the steady-state speeds $v_n$ we have: $v_n \le v_{n+1} \le v_{**}$ for all $n$, and $\lim_n v_n = v_{**}$.
\end{thm}
 
Proof of Theorem~\ref{thm-speed} is in Section~\ref{sec-speed-limit-proof}.

\section{Traveling wave shapes in the case of exponentially distributed independent jumps}
\label{sec-tws-exp}

\begin{thm} 
\label{th-exp-wave-existence}
Consider the special case $J(x) =1-e^{-x}$ and recall notation $v_* = (1+\sqrt{\lambda})^2$. Assume $\lambda>0$. Then:

(i) A TWS $\phi$ for the original system exists if and only if $v \ge v_*$. If $v>v_*$, then
\beql{eq-wave-exp-tail}
- \lim_{x\to\infty} [\log (1-\phi_x)]/x = \zeta(v),
\eeql
where $\zeta(v)$ is defined in \eqn{eq-v-gamma-exp-inverse}.

(ii) For any $v > v_*$ 
there exists a TWS $\phi$ for the system with moving left boundary at speed $v$, and it is such that \eqn{eq-wave-exp-tail} holds.

(iii) For any $v < v_*$ there exists a unique TWS $\phi$ for the system with moving right boundary at speed $v$.
\end{thm}

The proof of Theorem~\ref{th-exp-wave-existence} is in Section~\ref{sec-exp-proofs}.

\section{MFL speed results}
\label{sec-mfl-speed-results}

The following is a natural and standard notion of the average speed of a given MFL. 
We say that the average speed of an MFL $f(t)$, with a given initial state $f(0) \in \cm$, is lower (upper) bounded by $v$, if for any quantile $\nu \in (0,1)$ of $f(t)$ its average speed of progress is lower (upper) bounded by $v$, that is
$$
\liminf_{t\to\infty}~ (\limsup_{t\to\infty}) ~ \frac{1}{t}[q_\nu(f(t)) - q_\nu(f(0))] \le ~(\ge)~ v.
$$
An equivalent definition of a  lower (upper) speed bound $v$ is that for any $\epsilon >0$,
$$
\lim_{t\to\infty} f_{(v-\epsilon) t}(t) =0 ~~ \left(  \lim_{t\to\infty} f_{(v+\epsilon) t}(t) =1  \right).
$$
Speed $v$ is the average speed of an MFL, if it is both lower and upper bound.

Clearly an MFL average speed depends on the initial state, and can be arbitrarily large. (This follows from the existence
of a traveling wave with $\lambda=0$ and arbitrarily large speed $v$.) 

As a corollary of Theorem~\ref{prop-mckean-type} and Proposition~\ref{prop-speed-brw} we obtain the following

\begin{prop}
\label{prop-bmfl-speed}
The average speed of the BMFL is $v_{**}$.
\end{prop}

Proposition~\ref{prop-bmfl-speed} in turn is used to prove (in Section~\ref{sec-mfl-speed-lower-proof}) the following 

\begin{thm}
\label{thm-mfl-speed-lower}
The average speed of any MFL is lower bounded by $v_{**}$.
\end{thm}

In cases when the existence of traveling waves, including traveling waves with moving boundaries, can be established -- as we did in the case of exponential jump sizes -- those traveling waves can be used as lower and/or upper bounds for MFLs to obtain more precise MFL speed bounds. Specifically, the following fact can be used, which is a corollary of Lemma~\ref{lem-mfl-monotone}.

\begin{prop} 
\label{prop-bounds-by-waves}
Consider an MFL $f(\cdot)$.  
If there exists a traveling wave $\phi(\cdot)$ with speed $v$, and possibly with a lower [resp., upper] boundary with speed $v$,
such that $f(0) \preceq \phi(0)$ [resp., $\phi(0) \preceq f(0)$], then $f(t) \preceq \phi(t)$ [resp., $\phi(t) \preceq f(t)$]
at all times and, consequently, 
the MFL speed is lower [resp., upper] bounded by $v$.
\end{prop}

\begin{thm}
\label{thm-mfl-speed-exp}
Assume $J(x) = 1-e^{-x}$. 

(i) The average speed of the BMFL is $v_*$, and the average speed of any MFL is lower bounded by $v_{*}$.

(ii) If the right tail exponent of an MFL $f(\cdot)$ initial state $f(0)$ is lower bounded by $\zeta \le \zeta_*$, \\ i.e.
$\liminf_{x\to\infty} - \log (1-f_x(0))/x \ge \zeta$, then the MFL average speed is upper bounded by $v(\zeta)$. \\
Consequently, if $\liminf_{x\to\infty} - \log (1-f_x(0))/x \ge \zeta_*$, then 
the MFL average speed is $v_*$.

(iii) If the right tail exponent of an MFL $f(\cdot)$ initial state $f(0)$ is upper bounded by $\zeta \le \zeta_*$, \\ i.e.
$\limsup_{x\to\infty} - \log (1-f_x(0))/x \le \zeta$, then the MFL average speed is lower bounded by $v(\zeta)$.
\end{thm}

Proof of Theorem~\ref{thm-mfl-speed-exp} is in Section~\ref{sec-proof-mfl-speed-exp}. 
This proof does {\em not} rely on Proposition~\ref{prop-bmfl-speed}, which in turn relies on the connection to and properties of the associated BRW, described in Theorem~\ref{prop-mckean-type} and Proposition~\ref{prop-speed-brw}. 
It only relies on Theorem~\ref{th-exp-wave-existence} and Proposition~\ref{prop-bounds-by-waves}.
(This is the reason why we include statement (i) into Theorem~\ref{thm-mfl-speed-exp},
even though it is a special case of Proposition~\ref{prop-bmfl-speed} and Theorem~\ref{thm-mfl-speed-lower}.
We give an alternative proof of (i).) 

Therefore, the proof of Theorem~\ref{thm-mfl-speed-exp} shows how MFL speed bounds can be obtained directly from the existence and properties of traveling waves, which in our case are given by Theorem~\ref{th-exp-wave-existence}. However, the proof is quite generic. For example, some of the MFL speed bounds for the KPP model under certain initial conditions, obtained in \cite{larson-1978,freidlin-1979,freidlin-1985,veret-1999} via direct analysis or via Feynmann-Kac representation and large deviations analysis, can be instead obtained essentially the same way as our Theorem~\ref{thm-mfl-speed-exp}. Indeed, given the traveling waves' existence results for the KPP model \cite{kpp-1937},
 it is straightforward to obtain an analog of our Theorem~\ref{th-exp-wave-existence},
namely the existence of traveling waves with moving left and right boundaries for the speeds larger and smaller, respectively, than the critical speed. Then an analog of our Theorem~\ref{thm-mfl-speed-exp} follows, with essentially same proof.

\section{Conjectures, simulation results, discussion}
\label{sec-conclusions}

\subsection{Conjectures.}

We have proved that that the limit of the steady-state speeds $\lim_{n\to\infty} v_n = v_{**}$ holds for our model, under the general distribution of an independent jump size. However, we proved Theorem~\ref{th-exp-wave-existence} and Theorem~\ref{thm-mfl-speed-exp}(ii)-(iii) only for the exponential jumps.
We conjecture that analogs of Theorem~\ref{th-exp-wave-existence} and Theorem~\ref{thm-mfl-speed-exp}(ii)-(iii), in fact, hold for generally distributed jump sizes as well. We also conjecture that stationary distributions of centered processes converge to the Dirac measure concentrated on a TWS. Formally, we put forward the following

\begin{conjecture}
\label{conj-min-speed}
Suppose the jump size distribution $J(\cdot)$ is such that $\alpha >0$. Then:

(i) The unique TWS exists for any speed $v \ge v_{**}$, and does not exist for speeds $v< v_{**}$.

(ii) If the right tail exponent of an MFL $f(\cdot)$ initial state $f(0)$ is lower bounded by $\zeta \le \zeta_{**}$, \\ i.e.
$\liminf_{x\to\infty} - \log (1-f_x(0))/x \ge \zeta$, then the MFL average speed is upper bounded by $v(\zeta)$. \\
Consequently, if $\liminf_{x\to\infty} - \log (1-f_x(0))/x \ge \zeta_{**}$, then 
the MFL average speed is $v_{**}$.

(iii) If the right tail exponent of an MFL $f(\cdot)$ initial state $f(0)$ is upper bounded by $\zeta \le \zeta_{**}$, \\ i.e.
$\limsup_{x\to\infty} - \log (1-f_x(0))/x \le \zeta$, then the MFL average speed is lower bounded by $v(\zeta)$.

(iv) $\mathring f^n(\infty) \Rightarrow \phi^{**}$, where $\phi^{**}$ is the TWS for the speed $v_{**}$. 
\end{conjecture}

Proving Conjecture~\ref{conj-min-speed}, or at least some parts of it, is an interesting subject of future work.

\subsection{Simulations and discussion.}

Before presenting simulation results, we note that
if the system objective is to maximize the steady-state speed progress of a large-scale system, 
our Theorem~\ref{thm-speed} allows one to easily optimize the trade-off between $\lambda$ and $\mu$,
if they can be chosen subject to some constraint(s). For example, one may want to optimize the trade-off between efforts allocated to self-propulsion and synchronization, if both efforts require the consumption of a common limited resource. Our simulations show good accuracy of the optimal setting based on Theorem~\ref{thm-speed}.

Tables~\ref{sim-results-exp} and \ref{sim-results-unif} show steady-state speeds $v_n$, obtained by simulation of a system with $n=10000$ particles, under the exponentially and uniformly (in [0,2]) distributed jump sizes, respectively.
(Each simulation run is such that there are $400 \times n$ attempted jumps in total, i.e. $400$ per particle. This corresponds to simulation time
$400 \times n / (\lambda+\mu)$. First half of the simulation run is warmup.)
The tables also show the values of $v_{**}$. (For the exponential jump sizes, $v_{**}= v_*=(\sqrt{\lambda} + \sqrt{\mu})^2$.) 

The results for the exponential jump sizes are in Table~\ref{sim-results-exp}. We see that
the actual steady-state speed $v_n$ stays below $v_{**}$, as we know it should. We also observe that, even for a relatively large number of 10000 particles, there is still
a significant difference $v_{10000} - v_{**}$. This is not unusual, because the convergence $v_n \uparrow v_{**}$ 
in other contexts is known to be rather slow \cite{brunet-derrida-1997,brunet-derrida-2001}. Formally showing that 
this is true for our system as well is an interesting question for further research.
The results for the uniform jump sizes in Table~\ref{sim-results-unif} show the same patterns.

Suppose now that independent jumps and synchronizations by a particle consume some common resource (say, computing power or energy). To be specific, suppose the maximum rate at which the resource may be consumed by a particle is normalized to $1$, and the amounts of the resource consumed by one independent jump and one synchronization are $a>0$ and $b>0$, respectively. Suppose, it is desirable to maximize the speed of the particle system. 
Assume the number of particles is large so that $v_n \approx v_{**}$. 
Then we obtain the optimization problem 
\beql{eq-opt-prob}
\max v_{**} ~~\mbox{subject to} ~~ a \lambda + b \mu = 1.
\eeql
In the case of exponentially distributed jumps
$v_{**}= (\sqrt{\lambda} + \sqrt{\mu})^2$, and problem \eqn{eq-opt-prob}
has a simple explicit solution 
$$
\lambda_{opt} = \frac{1}{a+a^2/b}, ~~\mu_{opt} = \frac{1}{b+b^2/a}.
$$
For a general jump size distribution the optimal solution $(\lambda_{opt},\mu_{opt})$ can be easily found numerically using  \eqn{eq-v-gamma} and \eqn{eq-gamma-opt}.

In Tables~\ref{sim-results-exp} and \ref{sim-results-unif}, for the cases of exponential and uniform jumps, 
we vary $(\lambda,\mu)$, while satisfying constraint
$a\lambda+b\mu=1$ with $a=2, b=1$. In both cases we show $(\lambda_{opt},\mu_{opt})$. 
We see that, even though the values of $v_n$ for $n=10000$ are not necessarily very close to their limiting values $v_{**}$ ``yet'' (because this $n$ is not large enough), the parameter setting $(\lambda_{opt},\mu_{opt})$, obtained from \eqn{eq-opt-prob}, essentially maximizes $v_n$ for this $n$. Thus, the setting $(\lambda_{opt},\mu_{opt})$, which is easily computable, gives a practical ``rule of thumb'' for optimizing the trade-off between ``self-propulsion and synchronization'' in large-scale systems.

\begin{table}
\centering	
\begin{tabular}[]{| c | c || c | c |}
	\hline
	 $\lambda$ & $\mu$ & $v_n $ (simulation) & $v_{**}= (\sqrt{\lambda} + \sqrt{\mu})^2$ \\
	\hline
	$0.45$ & $0.1$ & $0.9321$ & $0.974264069$  \\
	$0.4$ & $0.2$ & $1.0863$ & $1.165685425$  \\
	$0.35$ & $0.3$ & $1.2104$ & $1.29807407$  \\
	$0.3$ & $0.4$ & $1.2974$ & $1.392820323$  \\
	$0.25$ & $0.5$ & $1.3236$ & $1.457106781$  \\
	$0.2$ & $0.6$ & $1.3318$ & $1.492820323$  \\
	$\lambda_{opt}=1/6$ & $\mu_{opt}=2/3$ & $1.3566$ & $1.5$  \\
	$0.15$ & $0.7$ & $1.3071$ & $1.49807407$  \\
	$0.1$ & $0.8$ & $1.2206$ & $1.465685425$  \\
	$0.05$ & $0.9$ & $1.0567$ & $1.374264069$  \\
	\hline
\end{tabular}
\caption{Exponential jump sizes, $n=10000$, $2\lambda+\mu=1$. Steady-state speed $v_n$ (simulated) and $v_{**}=\lim_n v_n$. (For exponential jump sizes, $v_{**} =  (\sqrt{\lambda} + \sqrt{\mu})^2$.) Pair $(\lambda_{opt}, \mu_{opt})$ maximizes $v_{**}$.} 
\label{sim-results-exp}
\end{table}

\begin{table}
\centering	
\begin{tabular}[]{| c | c || c | c |}
	\hline
	 $\lambda$ & $\mu$ & $v_n $ (simulation) & $v_{**}$ \\
	\hline
	$0.45$ & $0.1$ & $0.8176$ & $0.844$  \\
	$0.4$ & $0.2$ & $0.9243$ & $0.955$  \\
	$0.35$ & $0.3$ & $0.9704$ & $1.0165$  \\
	$0.3$ & $0.4$ & $0.9871$ & $1.0458$  \\
	$\lambda_{opt} \approx 0.27$ & $\mu_{opt} \approx 0.46$ & $0.9917$ & $1.0505$  \\
	$0.25$ & $0.5$ & $0.9995$ & $1.0486$  \\
	$0.2$ & $0.6$ & $0.9716$ & $1.0262$  \\
	$0.15$ & $0.7$ & $0.919$ & $0.9761$  \\
	$0.1$ & $0.8$ & $0.8209$ & $0.8907$  \\
	$0.05$ & $0.9$ & $0.6751$ & $0.7469$  \\
	\hline
\end{tabular}
\caption{Uniform$[0,2]$ jump sizes,  $n=10000$, $2\lambda+\mu=1$. Steady-state speed $v_n$ (simulated) and $v_{**}=\lim_n v_n$. Pair $(\lambda_{opt}, \mu_{opt})$ maximizes $v_{**}$.} 
\label{sim-results-unif}
\end{table}

\section{Proof of Theorem~\ref{th-finite-interval-B}}
\label{sec-MFL-B}

\subsection{Characterization of distributional limit.}
\label{sec-weak-limit}

\begin{thm}
\label{th-finite-interval-weak-limit}
Suppose $f^n(0) \stackrel{w}{\rightarrow} f(0)$, as in Theorem~\ref{th-finite-interval-B}.
Then the sequence of processes $f^n(\cdot)$ is tight in space $D([0,\infty), \cm^B)$, and any subsequential distributional limit $f(\cdot)$ is such that, w.p.1., the trajectory $f(\cdot)$ 
is a continuous element of $D([0,\infty), \cm^B)$,  
which satisfies \eqn{eq-pde-oper-B}.
\end{thm}

{\em Proof.} 
The proof is easily obtained by using the steps given in the proof of theorem 3 in \cite{St2022-wave-lst} (for a different particle system), namely theorems 4-7 and corollary 8 in \cite{St2022-wave-lst}. (In turn, the proof of theorem 3 in \cite{St2022-wave-lst} closely follows the development in \cite{balazs-racz-toth-2014}.) The only substantial difference is verification of condition (i) of \cite{perkins-2002}, which, in our notation, has the following form: for any $T>0$ and $\epsilon>0$, there exists  $K>0$ such that
\beql{eq-perkins-i}
\sup_n \pr \left(\sup_{t\le T}    [f^n_{-K}(t) + 1- f^n_{K}(t)] >\epsilon  \right) < \epsilon.
\eeql
However, for the process under consideration, with fixed right boundary at point $B$,
\eqn{eq-perkins-i} is trivial -- it suffices to choose $K$ large enough so that $K > B$
and  $f^n_{-K}(0) \le \epsilon/2$ uniformly in $n$. 
$\Box$

\subsection{Equivalent characterization of solutions to \eqn{eq-pde-oper-B} as mean-field models.}
\label{sec-mfm}

\begin{thm}
\label{lem-fde-oper-equiv}
For any initial condition $f(0) \in \cm^B$, a 
trajectory $f(\cdot) \in D([0,\infty), \cm^B)$ satisfies \eqn{eq-pde-oper-B} if and only if it is a mean-field model. 
\end{thm}

{\em Proof.} This proof is very similar to the proof of theorem 10 in \cite{St2022-wave-lst}. 

`Only if.' We only need to consider points $x<B$. Let $h=(h(u), ~u\in \R)$ be the 
step-function $h(u) = \bI(u \le x)$,
jumping from 1 to 0 at point $x$.
Let $h_\epsilon$, $\epsilon>0$, be a continuous approximation of $h$, which 
is linearly decreasing from $1$ to $0$ in $[x,x+\epsilon]$. We see that 
$$
L^B [f(t) h_\epsilon] \to L^B [f(t) h], ~~\forall t.
$$
Indeed, $|L^B [f(t) h_\epsilon] - L^B [f(t) h]|$ is upper bounded by $(\lambda+1)$ times the probability that a random jump -- either independent of size $Z$ or due to synchronization -- of a particle 
 randomly located according to distribution $f(t)$, is such that the jump either originates or lands in $(x,x+\epsilon)$;
 this probability vanishes as $\epsilon\downarrow 0$. Also, since both $L^B [f(t) h_\epsilon]$ and $L^B [f(t) h_\epsilon]$
 are within $[-(\lambda+1),0]$, for all $t$ and $\epsilon$, we have a universal bound $|L [f(t) h_\epsilon] - L [f(t) h]| \le \lambda+1$. 
 Then, for any fixed $t$,
by taking the $\epsilon\downarrow 0$ limit in $f(t) h_\epsilon - f(0) h_\epsilon - \int_0^t L^B [f(s) h_\epsilon] ds =0$, we obtain
$$
f(t) h  - f(0) h  - \int_0^t L^B [f(s) h] ds =0.
$$
This means that $f(t) h = f_x(t)$ is absolutely continuous in $t$, with the derivative equal to 
$(\partial/\partial t) f_x(t)=L^B [f(t) h]$ a.e. in $t$. This, in particular, implies that, for any fixed $y$, the $f_y(t)-f_{y-}(t)$
is continuous in $t$ (in fact, Lipschitz); this, in turn, means that, possible discontinuity points $y$ of $f_y(t)$
``cannot move''  in time $t$. We can then conclude that, for any $x$, the derivative
$(\partial/\partial t) f_x(t)=L^B [f(t) h]$ is in fact continuous in $t$.
Therefore, $(\partial/\partial t) f_x(t)=L^B [f(t) h]$ at every $t$. It remains to observe that
$L^B [f(t) h]$ is exactly the RHS of \eqn{eq-dyn}.

`If.'  The definition of a mean-field model $f(\cdot)$ implies that \eqn{eq-pde-oper-B} holds for the defined above step-function $h$ for any $x<B$.
Then, we have \eqn{eq-pde-oper-B} for any $h$,
which is piece-wise constant with finite number of pieces; and the set of such functions $h$ is tight within the space
of test functions $h \in \cc_b$, equipped with uniform metric. Then \eqn{eq-pde-oper-B} holds for any $h \in \cc_b$.
$\Box$

\subsection{MFM uniqueness.}
\label{sec-mfm-unique}

\begin{thm}
\label{lem-mfm-unique}
For any initial condition $f(0) \in \cm^B$,
an MFM  
$f(\cdot) \in D([0,\infty), \cm^B)$ is unique.
\end{thm}

{\em Proof.} The existence of an MFM for a given initial state $f(0)$ follows from
Theorems~\ref{th-finite-interval-weak-limit} and \ref{lem-fde-oper-equiv}.
For any MFM $f(\cdot)$ we know that $f_x(t)$, as a function of $t$, is Lipschitz, uniformly in $x$.
Consider two MFM trajectories, $f(\cdot)$ and $g(\cdot)$, with the same initial state $f(0)=g(0)$.
It is easy to check that for any $x$ and any $t>0$
\beql{eq-deriv-universal}
\left| \frac{d}{d t} g_x(t) - \frac{d}{d t} f_x(t) \right| \le C \| g(t) - f(t) \|
\eeql
for a fixed constant $C>0$ (for example, for $C=2\lambda+1$).
From \eqn{eq-deriv-universal} it is easy to see that 
\beql{eq-deriv-universal-2}
 \frac{d}{d t} \left\| g(t) - f(t) \right\| \le C \| g(t) - f(t) \|,
\eeql
as long as $\| g(t) - f(t) \| >0$. But then, by Gronwall inequality, $\| g(t) - f(t)\| \equiv 0$.
This proves the uniqueness.
$\Box$

\subsection{The conclusion of the proof of Theorem~\ref{th-finite-interval-B}.}
\label{sec-conv-th-proof-conclusion}

Theorem~\ref{th-finite-interval-B} follows immediately from the results of Subsections~\ref{sec-weak-limit} - \ref{sec-mfm-unique}.
$\Box$

\section{Proof of Theorem~\ref{th-finite-interval}}
\label{sec-MFL-orig}

We now generalize Theorem~\ref{th-finite-interval-B} to Theorem~\ref{th-finite-interval}. 
Recall that the evolution of $f_x^n(\cdot), x \in (-\infty, B),$ is independent of the behavior of particles located in $[B, \infty)$.
Theorem~\ref{th-finite-interval-B} holds for every $B$ for the projection $f^{n,B}(\cdot) \in D([0,\infty), \cm^B)$
of our original process $f^{n}(\cdot)$,
and we have convergence to the corresponding unique continuous 
MFL (MFM) $f^B(\cdot)$. By uniqueness, these MFLs for different $B$ must be consistent, namely for any $B_1 < B_2$, any $x< B_1$ and any $t$ we have $f_x^{B_1}(t)=f_x^{B_2}(t)$. So, we can formally define a trajectory
$f(t)$ by $f_x(t) \doteq \lim_{B\uparrow \infty} f_x^B(t)$. Let us show that this trajectory is the MFL satisfying
the statement of Theorem~\ref{th-finite-interval}. 

We first show that $f(t) \in \cm$ for any $t$, i.e. $f_\infty(t) = 1$. 
Recall that, for any $x$, $(d/dt)f_x(t)$ is non-positive, bounded and Lipschitz continuous in $t$ (uniformly in $x$).
Then $\epsilon(t) \doteq 1-f_\infty(t)$ is non-decreasing continuous, starting from $0$ at time $0$.
Moreover, $(d/dt)\epsilon(t)$ is bounded and Lipschitz continuous. 
Consider a sequence $c \uparrow \infty$ and the corresponding (non-decreasing) sequence of space-rescaled versions of $f(t)$:
$$
f^{(c)}_x(t) = f_{cx}(t).
$$
Clearly, $f^{(c)}(t)$ is the limiting trajectory of the process with rescaled initial state and rescaled jumps;
and for any $t$ the limit of $ f_{cx}(t)$, as a function of $x$, is constant, equal to $1-\epsilon(t)$.
Consider any fixed $x>0$  and $t>0$. As $c\to\infty$, $1- f^{(c)}_x(t) \to \epsilon(t)$ 
and $(d/dt)f^{(c)}_x(t) \to -\epsilon(t)(1-\epsilon(t))$. For fixed $t < s$, 
$$
\epsilon(s) - \epsilon(t) = \lim_{c\to\infty} [-f^{(c)}_x(s) + f^{(c)}_x(t)] 
= - \lim_{c\to\infty} \int_t^s (d/d\xi)f^{(c)}_x(\xi) d\xi = \int_t^s \epsilon(\xi)(1-\epsilon(\xi)) d\xi.
$$
We conclude that $\epsilon(t)$, as a function of $t$,
must be a solution to the logistic equation $\epsilon' = \epsilon(1-\epsilon)$ with initial condition $\epsilon(0)=0$.
Therefore, $\epsilon(t) \equiv 0$, and then $f_\infty(t)\equiv 0$, i.e. $f(t) \in \cm$. 

Now, because we know that $f(t) \in \cm$ for any $t$, $f(\cdot)$ is clearly an MFM -- 
it satisfies conditions of Definition~\ref{def-mfm} for any $x$ (because we can always choose $B >x$).
Condition \eqn{eq-pde-oper} holds for any $h \in \cc_b$, because we can choose $B$ large enough so that
$h(x)$ is constant in $[B,\infty)$ and for such $h$ the operators $L f(t) h = L^B f(t) h$, and then
\eqn{eq-pde-oper} follows from \eqn{eq-pde-oper-B}.

It remains to show the convergence
$f^n(\cdot) \Rightarrow f(\cdot)$. Fix any $T>0$, any $\epsilon>0$, and a sufficiently large $B=B(T)>0$, such that $1-f_B(T) \le \epsilon$. For the modified system with the right boundary at $B$, we have convergence of 
$f^{n,B}(\cdot) \in D([0,T], \cm^B)$ to a continuous $f^{B}(\cdot) \in D([0,T], \cm^B)$. But since $\epsilon>0$ can be arbitrarily small, we also have convergence of the original process 
$f^{n}(\cdot) \in D([0,T], \cm)$ to a continuous $f(\cdot) \in D([0,T], \cm)$. This is true for arbitrary $T>0$. Therefore, 
we also have convergence of the 
original process 
$f^{n}(\cdot) \in D([0,\infty), \cm)$ to a continuous $f(\cdot) \in D([0,\infty), \cm)$.

\section{Proof of Theorem~\ref{th-finite-interval-drift-1}}
\label{sec-MFL-1}

The proof is essentially same as that of Theorems~\ref{th-finite-interval} and \ref{th-finite-interval-B}, with straightforward adjustments. Condition \eqn{eq-perkins-i} for this process is automatic. 

We only comment on the proof of MFM uniqueness. Note that for any MFM $f(\cdot)$, for each $x$, the derivative $(d/dt) f_x(t)$
can have a discontinuity at time $t$ such that $B=x$ (i.e. when the moving boundary ``hits'' point $x$).
However, this does not affect the uniqueness proof in Section~\ref{sec-mfm-unique} -- we just need to note that
\eqn{eq-deriv-universal} and \eqn{eq-deriv-universal-2} hold for almost any $t>0$ w.r.t. Lebesgue measure. $\Box$

\section{Proof of Theorem~\ref{th-finite-interval-drift-2}}
\label{sec-MFL-2}

The proof is also essentially same as that of Theorems~\ref{th-finite-interval} and \ref{th-finite-interval-B}, with straightforward adjustments.  

Again, we only comment on the proof of MFM uniqueness. MFM for this system is such that, for a given $x$, the derivative $(d/dt) f_x(t)$ may have a discontinuity (jump down) at any time $t$ such that the jump distribution has an atom at $x-A$. 
If we consider two MFM trajectories, $f(\cdot)$ and $g(\cdot)$, with the same initial state $f(0)=g(0)$, 
as in the uniqueness proof in Section~\ref{sec-mfm-unique}, for each $x$, both $(d/dt) f_x(t)$ and $(d/dt) g_x(t)$
will have the same countable set of potential discontinuity time points. We also have that, at any $t$, $f_x(t)$ and
$g_x(t)$, as functions of $x$, have the same set of discontinuity points. Given these facts, we see that
the uniqueness proof in Section~\ref{sec-mfm-unique} applies -- again, 
we just need to note that
\eqn{eq-deriv-universal} and \eqn{eq-deriv-universal-2} hold for almost any $t>0$ w.r.t. Lebesgue measure. $\Box$

\section{Proof of Theorem~\ref{prop-mckean-type}}
\label{sec-mckean}

For the single particle independent jump (Markov) process $W(\cdot)$, 
denote by $P^t(x,H)$ its transition function and by $P^t$ its transition operator:
$$
P^t g_x = \int_x^\infty P^t(x,dy) g_y.
$$
The generator of this process is
\beql{eq-gener-W}
B g_x = \lambda \left[ \int_0^\infty  g_{x+\eta} d\bar J(\eta) - g_x \right].
\eeql
Any function $g$ such that $1-g \in \cm$ (i.e. a complementary distribution function) is within the generator $B$ domain,
and for such functions the integration by parts in \eqn{eq-gener-W} gives
\beql{eq-gener-W2}
B g_x = \lambda \int_x^\infty \bar J(\eta-x) d g_\eta.
\eeql
(The integration  in \eqn{eq-gener-W2} includes a possible atom of measure $g_\eta$ at $\eta=x$.
If $g$ is continuous, this subtlety is irrelevant.) 

Let us fix a function $g$ s.t. $1-g \in \cm$ and consider function
\beql{eq-branch-char}
u_x(t) = \E [g_{x+W_1(t)} \cdot \ldots \cdot g_{x+W_{N(t)}(t)}].
\eeql
Clearly, $u_x(0)=g_x$. Note that if 
$g_x = \bI\{x \le 0\} = \chi_{-x}$, then $u_x(t) = \pr\{\max_i W_i(t) \le 0\}$. Considering two cases -- when the first split of the first particles does not and does occur in $[0,t]$ -- we obtain
$$
u_x(t) = e^{-t} P^t g_x + \int_0^t e^{-s} P^{s} u^2_x(t-s) ds = e^{-t} P^t g_x + \int_0^t e^{-(t-\xi)} P^{t-\xi} u^2_x(\xi) d\xi,
$$
$$
e^t u_x(t) = P^t g_x + \int_0^t e^\xi P^{t-\xi} u^2_x(\xi) d\xi.
$$
Differentiating in $t$, obtain
$$
e^t u_x(t) + e^t \frac{\partial}{\partial t}u_x(t) = 
B P^t g_x
+ e^t  P^0 u^2_x(t) 
+ \int_0^t e^\xi P^{t-\xi} B u^2_x(\xi) d\xi,
$$
then
$$
\frac{\partial}{\partial t}u_x(t) = B u_x(t) + u^2_x(t) - u_x(t),
$$
and finally, recalling \eqn{eq-gener-W2},
\beql{eq-u-char}
\frac{\partial}{\partial t}u_x(t) = \lambda \int_x^\infty \bar J(\eta-x) d_\eta u_\eta(t) + u^2_x(t) - u_x(t).
\eeql
(It is easy to see from the definition of $u_x(t)$ that, for any $t>0$, it is continuous in $x$. Therefore, the subtlety 
mentioned immediately after \eqn{eq-gener-W2} is irrelevant.)

Now, let $g_x = \bI\{x \le 0\}$. Then, recall, $u_x(t) = \pr\{\max_i W_i(t) \le 0\}$, and then
$$
u_{-x}(t) = \pr\{\max_i W_i(t) \le x\}. 
$$
Substituting $f_x(t)=u_{-x}(t)$ in \eqn{eq-u-char} we obtain exactly \eqn{eq-dyn} with $f(0)=\chi$. $\Box$

\section{Proof of Theorem~\ref{thm-speed}.}
\label{sec-speed-limit-proof}

The monotonicity, $v_n \le v_{n+1}$, has already been observed in \eqn{eq-vn-increasing}.

Let us show that $v_n \le v_{**}$ for all $n$. 
The proof of this fact is same as the proof of the analogous upper bound in theorem 2.3(4) in \cite{jonck-fpp-2020}, except for our process we need to use Proposition~\ref{prop-speed-brw}, 
instead of the analogous fact for the process in \cite{jonck-fpp-2020} (which follows directly from the results of \cite{mckean-1975-kpp-branching} and \cite{kpp-1937}). Here is the proof, for completeness. 
Consider the system with finite number of $n$ particles. Consider another, artificial, system starting with the same initial state with $n$ particles, but such that each initial particle generates its own independent associated BRW as defined above. Clearly, the two systems can be coupled so that the artificial system dominates the actual one in the sense of $\preceq_l$, and in particular the location of its the right-most particle $M(t)$ is always to the right of the right-most particle location $M_n(t)$ of the actual system. But, Proposition~\ref{prop-speed-brw} holds for each of the $n$ independent associated BRW. We obtain $\limsup_{t\to\infty} M_n(t)/t \le \limsup_{t\to\infty} M(t)/t = v_{**}$ w.p.1, which implies $v_n \le v_{**}$.

It remains to show that $\liminf_n v_n \ge v_{**}.$ Consider arbitrary initial state of the process with $n$ particles.
Suppose WLOG that the leading particle is initially at $0$, and let
$D_n(t)$ be the location of the leading particle at time $t$. Consider a modified process, such that all particles, except the leading one are initially placed at $-\infty$. By monotonicity, the location $\tilde D_n(t)$ of the leading particle in the modified system is stochastically dominated by $D_n(t)$. The modified process is such that, in particular,
at some points in time some particles located at $-\infty$ with jump forward to join particles in $[0,\infty)$; let us call the particles located in $[0,\infty)$ regular. Let us now fix a finite time interval $[0,T]$, and consider the limit of the process of regular particles in this time interval, as $n\to\infty$. It is easy to construct coupling such that, w.p.1, the process of regular particles (of modified system) converges to an associated BRW. (Indeed, when $n$ is large, the rate at which each regular particle is being joined by a particle from $-\infty$ is close to $1$. We omit details, which are straightforward.) By Proposition~\ref{prop-speed-brw}, if $T$ is sufficiently large, the average speed of advance 
of the leading particle of the BRW in $[0,T]$ is close to $v_{**}$. (This is if $\alpha>0$ and then $v_{**} < \infty$.
If $\alpha=0$ and then $v_{**} = \infty$, ``close to $v_{**}$'' means ``arbitrarily large.'')
Consequently, for all large $n$, the average speed of advance of $\tilde D_n(t)$ is close to $v_{**}$. 
This implies (we skip straightforward $\epsilon-\delta$ formalities)
that $\liminf_n v_n \ge v_{**}$. $\Box$.

\section{Existence and properties of traveling waves in the case of exponentially distributed jump sizes. Proof of Theorem~\ref{th-exp-wave-existence}} 
\label{sec-exp-proofs}

In this section we consider the case of exponential jump size distribution,
$J(x) =1- e^{-x}$, and develop results leading to the proof of Theorem~\ref{th-exp-wave-existence}.

Recall, that in this case a TWS $\phi=(\phi_x)$ must satisfy ODE
\eqn{eq-tws-exp-1}:
\beql{eq-tws-exp}
v \phi'' = (1+\lambda-v-2\phi) \phi'  + \phi(1-\phi).
\eeql
A function $\phi=(\phi_x)$ is a TWS for our original system 
if and only if it is a {\em proper} solution of ODE \eqn{eq-tws-exp}, namely: 
$\phi_x \in [0,1]$ for all $x \in \R$,
\eqn{eq-tws-exp} holds for all $x \in \R$,
$\lim_{x \downarrow -\infty} \phi_x = 0$ and $\lim_{x \uparrow \infty} \phi_x = 1$.

It will be convenient to consider ODE \eqn{eq-tws-exp} in terms of the first order ODE in phase space $(\phi, z=\phi')$:
\beql{eq-tws-exp-phase}
\phi' = z, ~~z' = -z + \frac{1+\lambda-2\phi}{v}z + \frac{\phi}{v}(1-\phi).
\eeql
In terms of \eqn{eq-tws-exp-phase}, a proper solution is a solution which starts at point $(0,0)$, ends 
at point $(1,0)$, and stays within the strip $\phi \in [0,1]$.

The following are immediate observations from \eqn{eq-tws-exp-phase}. Solutions to \eqn{eq-tws-exp-phase} (not necessarily proper solutions) are such that: a solution trajectory cannot hit $\phi$-axis for $\phi\in (0,1)$; 
any two solution trajectories either coincide in $\phi \in [0,1]$ or one strictly dominates another for $\phi \in (0,1)$.

\subsection{Speed $v$ lower bound, for any proper solution}

\begin{figure}
\centering
        \includegraphics[width=4in]{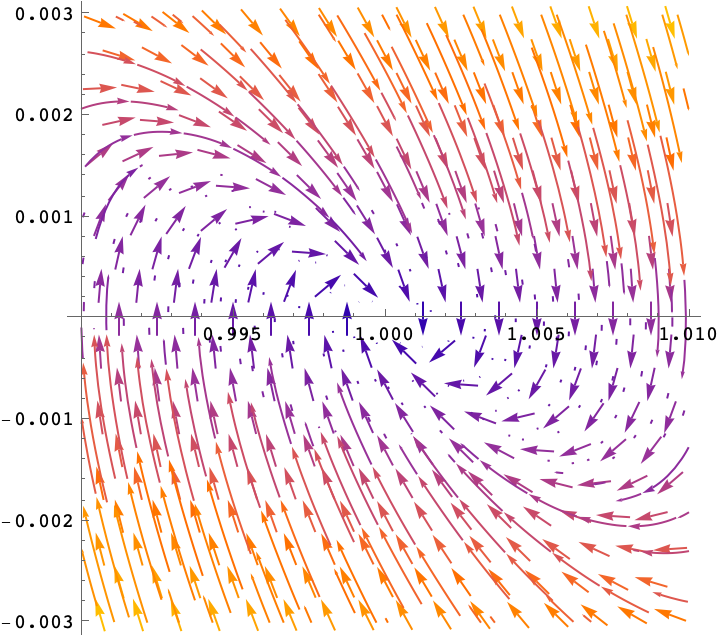}
    \caption{Vector-field $(\phi',z')$ in the vicinity of $(1,0)$; $\lambda=4$, $v=7 < v_*=9$.}
    \label{fig-image-10}
\end{figure}

Consider the linearization of \eqn{eq-tws-exp-phase} at point $(1,0)$. Denoting $\nu=1-\phi$, the linearization 
is:
\beql{eq-tws-exp-phase-linear-end}
\nu' = -z, ~~z' = -z + \frac{\lambda-1}{v}z + \frac{1}{v}\nu.
\eeql
The eigenvalues of this linear system satisfy characteristic equation
\beql{eq-tws-exp-tail-char}
v \zeta^2 + (1+v -\lambda) \zeta + 1 =0.
\eeql
We know (see Section~\ref{sec-trav-waves})
 that for a proper solution to exist, we must have $v\ge \lambda$. 
Given that, $-(1+v-\lambda) < 0$, and it is easy to check that condition $(1+v-\lambda)^2 -4v \ge 0$ 
is equivalent to condition $v \ge v_* \doteq (1+\sqrt{\lambda})^2$.
We see that: both eigenvalues of \eqn{eq-tws-exp-phase-linear-end} (roots of \eqn{eq-tws-exp-tail-char}) have negative real parts; the eigenvalues are two (different) adjoint complex numbers if and only if $v<(1+\sqrt{\lambda})^2$, and they are both real otherwise. We then observe that 
when $v<(1+\sqrt{\lambda})^2$, a proper solution does {\em not} exist, because any solution 
with end point at $(1,0)$ converges to $(1,0)$ by oscillating around it (as illustrated in Figure~\ref{fig-image-10}) and thus 
must ``hit'' the line $\phi=1$ at a point strictly above point $(1,0)$.

We conclude that: {\em if a proper solution of \eqn{eq-tws-exp-phase} exists, then necessarily}
$
v\ge v_* = (1+\sqrt{\lambda})^2.
$

\subsection{Proof of Theorem~\ref{th-exp-wave-existence}(i).}

To complete the proof of Theorem~\ref{th-exp-wave-existence}(i) we need to show that a 
a proper solution of \eqn{eq-tws-exp-phase} exists for any $v \ge v_*$. We will prove that it exists for 
$v> v_*$, and then will obtain the existence for $v=v_*$ by continuity.

\subsubsection{Existence/uniqueness/continuity of solutions to \eqn{eq-tws-exp-phase}.}

Here we consider solutions to \eqn{eq-tws-exp-phase} (not necessarily proper solutions) 
within the domain $\{0\le \phi \le 1, ~z\ge 0\}$. Consider a fixed parameter $v \ge v_*$ and an initial condition 
$(\phi_0,z_0)$ for a solution \eqn{eq-tws-exp-phase}. (Initial condition can be viewed as a point which a solution trajectory must contain, because the trajectories satisfying \eqn{eq-tws-exp-phase} can be considered both in the 
forward and reverse directions.) From the standard ODE theory, we see that the solution for a given triple
$(\phi_0,z_0,v)$ exists and is unique, as long as $(\phi_0,z_0)\ne (0,0)$ and $(\phi_0,z_0)\ne (1,0)$.
Moreover, the solution is continuous in $(\phi_0,z_0,v)$ outside those two points.
The following lemma shows that the existence/uniqueness/continuity holds for the initial condition
$(\phi_0,z_0)= (0,0)$ as well, if we consider non-trivial solutions.

\begin{lem}
\label{lem-uniquness-finite}
For a fixed $\lambda >0$, and any speed $v\ge v_*$, there exists a unique non-trivial solution  $(\phi, z(\phi))$, such that $z(\phi)>0$ and $(\phi, z(\phi)) \to (0,0)$ as $\phi\downarrow 0$.
This solution is such that $dz/d\phi = \gamma$ at initial $\phi=0$, where 
\beql{eq-gamma-def}
\gamma=\frac{-(v-1-\lambda) + \sqrt{(v-1-\lambda)^2 +4v}}{2v} >0.
\eeql
Moreover, this solution is such that the solution continuity w.r.t. $(\phi_0,z_0,v)$ holds.
\end{lem}

\begin{figure}
\centering
        \includegraphics[width=4in]{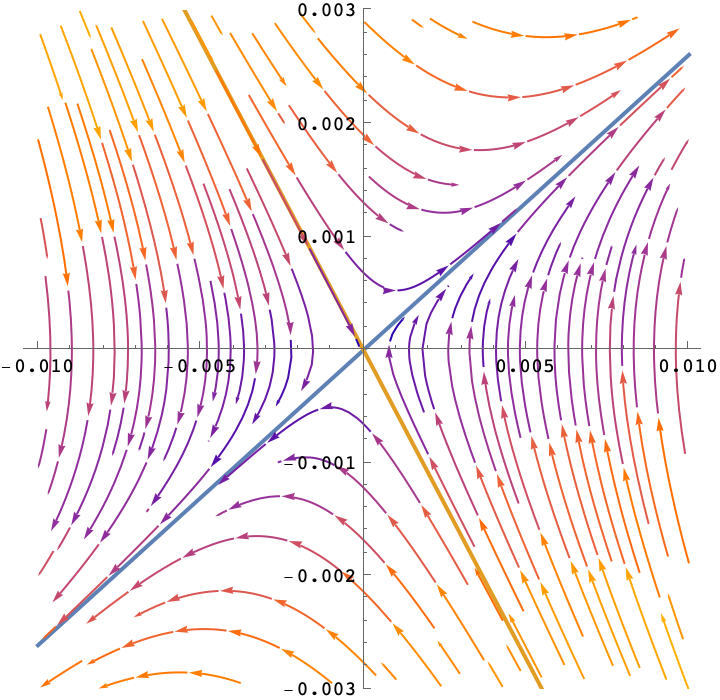}
    \caption{Vector-field $(\phi',z')$ in the vicinity of $(0,0)$; $\lambda=4$, $v=10$, $v > v_*$.}
    \label{fig-image-00}
\end{figure}

{\em Proof of Lemma~\ref{lem-uniquness-finite}.} 
Recall that $v \ge v_* > 1+ \lambda$. 
Then, for the linearization of \eqn{eq-tws-exp-phase} at point $(0,0)$, we always have two different real eigenvalues, one positive and one negative: $-\gamma_2 < 0 < \gamma$, where $\gamma$ is given by \eqn{eq-gamma-def}.
The corresponding eigenvectors are $(1,\gamma)$ and $(1,-\gamma_2)$, see Figure~\ref{fig-image-00}.

The proof of existence is constructive. Consider any sequence of solutions, starting at a points $(\phi^{(k)}_0, z^{(k)}_0) \ge (0,0)$ and corresponding to a speed $v^{(k)} \ge v_*$; the sequence is such that, as $n\to\infty$, $(\phi^{(k)}_0, z^{(k)}_0) \to (0,0)$, $v^{(k)} \to v$, and $(\phi^{(k)}_0, z^{(k)}_0) \ne (0,0)$ for all $k$. 
We view these solutions as functions $z^{(k)}=z^{(k)}(\phi)$; this is possible because $\phi'=z>0$ 
for $\phi \in (\phi_0, 1)$.
Then any sub-sequential limit of this sequence of solutions must be a solution $z=z(\phi)$ (with speed parameter $v$)
starting from $(0,0)$. Also, given the structure of the derivative vector-field around $(0,0)$, it is easy to see that any such solution must satisfy $dz/d\phi = \gamma$ at $\phi=0$. 
(Indeed, for $c>0$ denote by $M_c = \{z=c\phi, ~\phi \ge 0\}$ the ray from point $(0,0)$ in the direction $(1,c)$.
For a small fixed $\epsilon>0$, denote by $W$ the closed cone between the rays $M_{\gamma-\epsilon}$ and 
$M_{\gamma+\epsilon}$. If we look at the solution in the reverse direction, as it approaches point $(0,0)$,
it cannot be outside $W$ arbitrarily close to $(0,0)$ -- otherwise it will hit $z$-axis strictly above $(0,0)$.)

It remains to show the uniqueness of a solution for a given parameter $v$. (Then the continuity in $(\phi_0,z_0,v)$ will automatically follow.) Consider any two solutions to \eqn{eq-tws-exp-phase}, viewed as 
functions $z=z(\phi)$ and $\hat z= \hat z(\phi)$. It is easy to see from \eqn{eq-tws-exp-phase} that
$$
\frac{d}{d\phi} [\hat z - z] = \frac{\phi(1-\phi)}{v} \left[ \frac{1}{\hat z}  -   \frac{1}{z} \right]
= - \frac{\phi(1-\phi)}{v\hat z z}  [\hat z - z].
$$
Thus, the difference $\hat z(\phi) - z(\phi)$ cannot increase in $\phi$, hence uniqueness.
$\Box$

\subsubsection{Proper solution existence for $v \ge v_*$.}

\begin{lem}
\label{lem-v-le-vstar}
A proper solution of \eqn{eq-tws-exp-phase} exists for any $v>v_*$.
(And then, by continuity of solutions in $v$, it exists for $v=v_*$ as well.)
\end{lem}

\begin{figure}
\centering
        \includegraphics[width=4in]{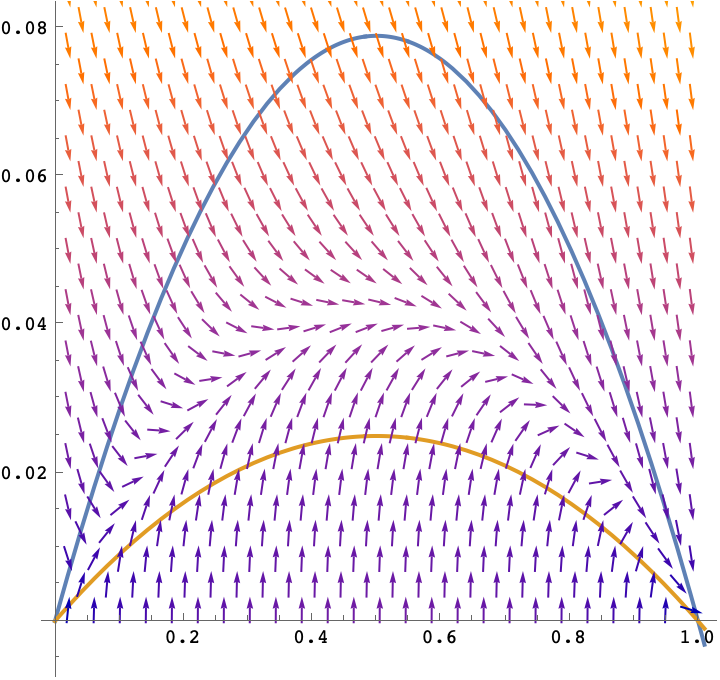}
    \caption{Vector-field $(\phi',z')$ points inside the domain between the parabolas $z = C \phi(1-\phi)/v$ and $z=\phi(1-\phi)/v$; here
    $C=\sqrt{v}$, $\lambda=4$, $v=10$, $v > v_*$.}
    \label{fig-image-parabolas}
\end{figure}

{\em Proof of Lemma~\ref{lem-v-le-vstar}.} 
Consider $C>1$ (to be determined later). Consider a point $(\phi,z)$, with $\phi\in (0,1)$ and $z\ge 0$.
Let us compare the derivative vector $(\phi',z')$, given by \eqn{eq-tws-exp-phase}, to the tangent vector to 
the parabola $h(\phi) = C \phi(1-\phi)/v$ at point $\phi$, with the $\phi$-component equal to $\phi'$;
the tangent vector is then $(\phi', \phi' C(1-2\phi)/v)$. The difference of $z$-components of these two vectors is
\[
z' - \phi' C(1-2\phi)/v =
-\left(\frac{v-\lambda+(C-1)(1-2\phi)}{v}-\frac{1}{C}\right)z-\frac{1}{C}(z-h(\phi)).
\]
If $(\phi,z)=(\phi, h(\phi))$, i.e. point $(\phi,z)$ is on the parabola, this difference is
$$
-\left(\frac{v-\lambda+(C-1)(1-2\phi)}{v}-\frac{1}{C}\right)z.
$$
Let us consider the expression in the bracket (a linear function of $\phi$). If we can prove that it is positive for $0\leq \phi\leq 1$, then {\em the trajectories of our dynamic system will traverse the parabola $h=C\phi(1-\phi)/v$ entering the domain between that parabola and $\phi$-axis.} (See Figure~\ref{fig-image-parabolas}.)
Let us check that this is the case.
The coefficient in front of $\phi$ is negative, so it suffices to check the expression is positive at $\phi=1$, i.e. 
\[
\frac{v-\lambda-(C-1)}{v}-\frac{1}{C}>0, \mbox{ i.e., }\frac{v-\lambda+1}{v}>\frac{C}{v}+\frac{1}{C}. 
\]
Condition $v>(1+\sqrt{\lambda})^2$ is equivalent to $\sqrt{v}>1+\sqrt{\lambda}$,
then to $(\sqrt{v}-1)^2>\lambda$, then to $v-\lambda+1>2\sqrt{v}$.
If we take $C=\sqrt{v}$, this gives exactly the last display. 

Now, consider the solution $z(\phi)$ to \eqn{eq-tws-exp-phase} with initial condition $(0,0)$. The derivative $dz/d\phi$
at $\phi=0$ is $\gamma$ (in \eqn{eq-gamma-def}), while the derivative of $h(\phi)$ (with $C=\sqrt{v}$) is $1/\sqrt{v}$.
For $\lambda>0$ and $v \ge v_*$, the inequality 
$1/\sqrt{v} > \gamma$
 is verified by simple algebra. This means that, for sufficiently small $\phi>0$, the solution $z(\phi)$ is ``under''
the parabola $h(\phi)$. But, as shown above, it cannot traverse $h(\phi)$. We also know that solution 
$z(\phi)$ cannot hit $\phi$-axis while $\phi\in (0,1)$. We conclude that $z(\phi) \le h(\phi)$ 
for all $\phi \le 1$, and $z(1)=0$. In particular, $z(\phi)$ is a proper solution.
(In fact, one can show that the trajectories also are transversal to the parabola $z=\phi(1-\phi)/v$.)
%$\Box$

\subsubsection{Proper solution derivative at point $(1,0)$.}

To complete the proof of of Theorem~\ref{th-exp-wave-existence}(i), it remains to show that \eqn{eq-wave-exp-tail}
holds when $v > v_*$. If $v>v_*$, the linearization of our dynamic system at point $(1,0)$ has two strictly negative eigenvalues, $-\zeta_1 > -\zeta_2,$, being roots of \eqn{eq-tws-exp-tail-char}, so that
$$
\zeta_1  = \frac{(1+v-\lambda) - \sqrt{(1+v-\lambda)^2 - 4v}}{2v}, ~~~ 
\zeta_2  = \frac{(1+v-\lambda) + \sqrt{(1+v-\lambda)^2 - 4v}}{2v}.
$$
Since solution $z=z(\phi)$ for speed $v$ is dominated by parabola $\phi(1-\phi)/\sqrt{v}$, then in the neighborhood of $(1,0)$ it is dominated by the line $z=(1-\phi)/\sqrt{v}$. It is easy to see that $1/\sqrt{v} < \zeta_2$. (Indeed, a stronger inequality $1/\sqrt{v} < (1+v-\lambda)/(2v)$ is equivalent to $v > v_*$.) Therefore, in the neighborhood of $(1,0)$
the solution $z=z(\phi)$ is separated from the eigendirection $z=\zeta_2 (1-\phi)$ by line $z=(1-\phi)/\sqrt{v}$.
But then the solution $z=z(\phi)$ must approach point $(1,0)$ along the eigendirection $z=\zeta_1 (1-\phi)$,
that is 
$$
\left. \frac{dz}{d\phi} \right|_{\phi=1} = - \zeta_1.
$$
It remains to notice that $\zeta_1$ is exactly equal to $\zeta(v)$ in \eqn{eq-v-gamma-exp-inverse}.
$\Box$

\subsection{Existence of solutions, giving TWS for systems with moving left or right boundary. Proof of Theorem~\ref{th-exp-wave-existence}(ii)-(iii).}

The following Theorem~\ref{thm-sol-left-boundary} is an extended version of Theorem~\ref{th-exp-wave-existence}(ii).

\begin{thm}
\label{thm-sol-left-boundary}
Suppose $\lambda>0$ and $v>v_*$.
For any sufficiently small $\phi_0>0$ and $z_0=\phi_0(1-\phi_0 + \lambda)/v$, the solution to \eqn{eq-tws-exp-phase} with initial condition $(\phi_0,z_0)$ ends at $(1,0)$.
This solution corresponds to the unique (up to a shift)  TWS $\phi=(\phi_x)$ for the system with moving left 
boundary at speed $v$, satisfying  conditions $\phi_{x_0}=\phi_0$, $\phi'_{x_0}=z_0$, where $x_0$ is the left boundary
point of the TWS; moreover, the TWS $\phi$ satisfies \eqn{eq-wave-exp-tail}.
\end{thm}

{\em Proof of Theorem~\ref{thm-sol-left-boundary}.}  Since $v>v_*$, there is the unique solution $z(\phi)$
to \eqn{eq-tws-exp-phase} with initial condition $(0,0)$, and it is proper. Moreover,
this solution is ``under'' the parabola $h(\phi)=\sqrt{v} \phi(1-\phi)/v$, as shown in the proof 
of Lemma~\ref{lem-v-le-vstar}. 
Let us consider the solution $\hat z(\phi)$ to \eqn{eq-tws-exp-phase} with initial condition $(0,\epsilon)$ with small $\epsilon>0$. By continuity of solutions to \eqn{eq-tws-exp-phase} in initial condition, 
we can choose $\epsilon$ small enough, 
so that $\hat z(\phi)$ is close to $z(\phi)$ and therefore is under the $h(\phi)$ for some values $\phi \in (0,1)$.
We conclude (as in the proof 
of Lemma~\ref{lem-v-le-vstar}) that solution $\hat z(\phi)$ for $\phi \in [0,1]$ ends at point $(1,0)$. We also know that 
$z(\phi) < \hat z(\phi)$ for all $\phi \in [0,1)$. Recall that the derivative $dz/d\phi$ at $\phi=0$ is equal to 
$\gamma$ in \eqn{eq-gamma-def}. 
The dependence $z_0=\phi_0(1-\phi_0 + \lambda)/v$ on $\phi_0$
has derivative $(1+\lambda)/v$ at $\phi_0=0$. The inequality 
$ 
(1+\lambda)/v > \gamma
$ 
is verified by simple algebra. We conclude the following: for all sufficiently small positive $\phi_0$,
the point $(z_0,\phi_0)$ lies strictly between the solutions $z(\phi)$ and $\hat z(\phi)$.
Therefore, the unique solution to \eqn{eq-tws-exp-phase} with initial condition $(z_0,\phi_0)$ ends at point $(1,0)$.

Verification of \eqn{eq-wave-exp-tail} for the corresponding TWS $\phi$ is the same as in the proof of
Theorem~\ref{th-exp-wave-existence}(i).
$\Box$

{\em Proof of Theorem~\ref{th-exp-wave-existence}(iii).} There is the unique solution to \eqn{eq-tws-exp-phase}
with initial condition $(0,0)$. Since $v<v_*$, this solution cannot be proper. Therefore, the solution 
``hits'' the line $\phi=1$ strictly above the $\phi$-axis, i.e. at a point $(1,z_1)$ with $z_1>0$.
The corresponding (unique up to a shift) distribution function $\phi=(\phi_x)$ is the unique TWS for the system 
with moving right boundary at speed $v$.
$\Box$

\section{Proof of Theorem~\ref{thm-mfl-speed-lower}.}
\label{sec-mfl-speed-lower-proof}

Consider a fixed MFL $f(\cdot)$. Let us fix an arbitrarily small $\nu>0$. WLOG assume that $f_0(0) \ge \nu$ and $f_{0-}(0) \le \nu$. Recall that trajectory $f(\cdot)$ is the unique distributional limit of a sequence of processes $f^n(\cdot)$
with deterministic initial states such that $f^n(0) \stackrel{w}{\rightarrow} f(0)$. We can and will choose a sequence of
initial states such that, for all $n$, $\nu n$ of left-most particles are located in $(-\infty,0]$ (let us call them ``left particles'') and then the remaining $(1-\nu)n$ particles are located in $[0,\infty)$ (let us call them ``right particles''). 

Now, for each $n$, consider a lower-bounding process $\hat f^n(\cdot)$ constructed as follows. 
The initial state of $\hat f^n(\cdot)$ is such that the $\nu n$ ``left'' particles' locations are the same as in the original process, i.e. $\hat f_x^n(0) = f_x^n(0)$ for $x< 0$, while the $(1-\nu) n$ ``right'' particles are located exactly at $0$.
The evolution of $\hat f^n(\cdot)$ is such that the ``left'' particles are ``frozen'' -- they do not jump at all;
the ``right'' particles jump as usual, making both independent and synchronization jumps. Clearly, 
$\hat f^n(\cdot)$ and $f^n(\cdot)$ can be coupled so that $\hat f^n(t) \preceq f^n(t)$ at all times; then 
$\hat f(t) \preceq f(t)$ for the corresponding MFLs.
Now, observe that  the evolution of $\hat f^n(\cdot)$ is such that the left particles can be simply ignored,
and all right particles are initially at $0$. Next observe that process $\hat f^n(\cdot)$ is such that the evolution of the set of $(1-\nu)n$ right particles can be viewed as that of a stand-alone original system, except the synchronization rate of each particle is $(1-\nu)$ instead of $1$. (Recall that in the system under consideration $\mu=1$, WLOG.)
Therefore, the MFL  $\hat f(\cdot)$ has the following form:
$\hat f_x(t) =f_x(0)(t)$ for $x<0$, and 
\beql{eq-tildef}
\hat f_x(t) =\nu + (1-\nu) \tilde f_x(t), ~~\mbox{for}~ x>0,
\eeql
where $\tilde f(\cdot)$ is the BMFL of the system with independent jump rate $\lambda$ and synchronization rate $1-\nu$. For a sufficiently small $\nu>0$, the speed $v'_{**}$ of $\tilde f(\cdot)$ is arbitrarily close to $v_{**}$. 
(In particular, if $\alpha=0$, then $v'_{**}=v_{**}=\infty$.)
We conclude
that any quantile $\beta > \nu$ of $f(\cdot)$ advances at  average speed at least $v'_{**}$. Since $\nu>0$ can be chosen 
arbitrarily small, any quantile $\beta > 0$ of $f(\cdot)$ advances at  average speed at least $v_{**}$. 
$\Box$

\section{Proof of Theorem~\ref{thm-mfl-speed-exp}.} 
\label{sec-proof-mfl-speed-exp}

(i) Obviously, this property is a special case of the more general properties, given in Proposition~\ref{prop-bmfl-speed}
and Theorem~\ref{thm-mfl-speed-lower}. However, we now give a simple alternative proof, which does {\em not}
rely on the connection to and properties of the associated BRW, described in Theorem~\ref{prop-mckean-type} and Proposition~\ref{prop-speed-brw}. The proof relies only on Theorem~\ref{th-exp-wave-existence} and monotonicity. 

By Theorem~\ref{th-exp-wave-existence}(ii)-(iii), the BMFL can be lower and upper dominated by the following traveling waves: a traveling wave with moving upper boundary with speed $v'$ slightly smaller that $v_*$,
and a traveling wave with moving lower boundary with speed $v''$ slightly larger that $v_*$. Since $v'$ and $v''$ can be arbitrarily close to $v_*$, the average speed of BMFL is $v_*$. To prove that $v_*$ is a lower bound of any MFL average speed, the proof of Theorem~\ref{thm-mfl-speed-lower} applies.

(ii) By Theorem~\ref{th-exp-wave-existence}(ii) there exists a TWS $\psi$ with moving left boundary and speed $v'$ slightly greater than $v(\zeta)$. The right tail exponent of $\psi$ is slightly smaller than $\zeta$. We can always consider a version of this $\psi$, shifted sufficiently far to the right, so that $f(0) \preceq \psi$. By monotonicity, MFL $f(\cdot)$ is dominated by the traveling wave (with moving left boundary) with speed $v'$, starting from $\psi$.
Therefore, the average speed of the MFL is upper bounded by $v'$, which can chosen arbitrarily close to $v(\zeta)$.

(iii) When $\zeta = \zeta_*$, the claimed property follows from (i). 
Therefore, it suffices to consider the case $\zeta < \zeta_*$. 
Fix a small $\nu>0$. If $\nu$ is small enough, by Theorem~\ref{th-exp-wave-existence}(iii) there exists a TWS $\psi$ for the system with synchronization rate $\mu=1-\nu$, with speed $v'$ slightly smaller than $v(\zeta)$, and with the right tail exponent of $\psi$ slightly larger than $\zeta$. Consider now function $\tilde \psi = \nu + (1-\nu) \psi$, and consider its version shifted sufficiently far to the left, so that $\tilde \psi \preceq f(0)$. This $\tilde \psi$ can be viewed as the initial state of an MFL for a system where $\nu n$ left-most particles are located at $-\infty$ and are frozen -- never make any jumps; then the corresponding MFL is a traveling wave with speed $v'$, dominated by $f(\cdot)$. We see that any quantile $\beta > \nu$ of $f(t)$ advances at  average speed at least $v'$. Since $\nu$ can be arbitrarily small and $v'$ can be arbitrarily close to $v(\zeta)$,
we obtain the claim.
$\Box$

%\iffalse
%%%%%%%%%%%%\bibliographystyle{acmtrans-ims}
%%%%%%%%%%%%\bibliographystyle{apt}
\bibliographystyle{abbrv}

%\bibliography{biblio-stolyar}
%\fi

\end{document}